\documentclass[a4paper, 10pt]{article}

\usepackage{amsmath,amsfonts,amssymb,amsthm,enumerate,graphicx}
\usepackage{amscd, amssymb}
\usepackage{tikz,authblk}
\usepackage{tikz-cd}
\usepackage{subfig}
\usepackage{url}
\usepackage{hyperref}
\usepackage{float}
\usepackage{xcolor}
\usepackage{marginnote}
\usepackage{easyReview}
\usepackage{boldline}
\hypersetup{
	colorlinks=true,
	citecolor=magenta,
	linkcolor=red,
	filecolor=black,      
	urlcolor=magenta}
	
\usepackage[T1]{fontenc}
\usepackage{ae}
\usepackage[all]{xy}

\usepackage{lscape}

\DeclareMathOperator{\Aut}{{\rm Aut}}
\DeclareMathOperator{\Nor}{{\rm Nor}}

\DeclareFontFamily{U}{BOONDOX-calo}{\skewchar\font=45 }
\DeclareFontShape{U}{BOONDOX-calo}{m}{n}{
  <-> s*[1.05] BOONDOX-r-calo}{}
\DeclareFontShape{U}{BOONDOX-calo}{b}{n}{
  <-> s*[1.05] BOONDOX-b-calo}{}
\DeclareMathAlphabet{\mathcalboondox}{U}{BOONDOX-calo}{m}{n}
\SetMathAlphabet{\mathcalboondox}{bold}{U}{BOONDOX-calo}{b}{n}
\DeclareMathAlphabet{\mathbcalboondox}{U}{BOONDOX-calo}{b}{n}
\usepackage[utf8]{inputenc}
         \DeclareMathAlphabet{\mathscr}{U}{BOONDOX-cal}{m}{n}
         \SetMathAlphabet{\mathscr}{bold}{U}{BOONDOX-cal}{b}{n}
         \DeclareMathAlphabet{\mathbscr} {U}{BOONDOX-cal}{b}{n}

\newtheorem{theorem}{{\bf Theorem}}[section]
\newtheorem{result}{{\bf Result}}[section]

\newtheorem{definition}[theorem]{{\bf Definition}}
\newtheorem{example}[theorem]{{\bf Example}}
\newtheorem{lemma}[theorem]{{\bf Lemma}}

\newtheorem{proposition}[theorem]{{\bf Proposition}}
\newtheorem{remark}[theorem]{{\bf Remark}}

\newtheorem{question}[theorem]{{\bf Question}}

\DeclareMathOperator{\Lim}{\limsup_{\textit{v} \to \infty}}

\newcommand{\lk}[2]{{\rm lk}_{#1}(#2)}
\newcommand{\skel}[2]{{\rm skel}_{#1}(#2)}
\newcommand{\Kd}[1]{{\mathcal K}(#1)}

\newcommand{\Kc}{\mathcal{K}}
\newcommand{\Hc}{\mathcal{H}}
\newcommand{\vg}{\vspace{0.5cm}}

\newcommand{\bs}{\backslash}
\newcommand{\ol}{\overline}
\newcommand{\meets}{\leftrightarrow}
\newcommand{\nmeets}{\nleftrightarrow}

\newcommand{\FF}{ \ensuremath{\mathbb{F}}}
\newcommand{\QQ}{ \ensuremath{\mathbb{Q}}}
\newcommand{\ZZ}{ \ensuremath{\mathbb{Z}}}
\newcommand{\RR}{ \ensuremath{\mathbb{R}}}

\newcommand{\IntA}{A^{\!\!\!^{^{\circ}}}}
\newcommand{\IntB}{B^{\!\!\!^{^{\circ}}}}

\newcommand{\TPSS}{S^{\hspace{.2mm}2}\! \times \hspace{-3.3mm}_{-} \,
S^{\hspace{.1mm}1}}

\renewcommand{\Authfont}{\scshape\small}
\renewcommand{\Affilfont}{\itshape\small}
\renewcommand{\Authand}{,}
\renewcommand{\Authands}{, }
\date{\today}

\begin{document}

\title{2-uniform toroidal maps, classification and asymptotic behavior}

	\author {Arnab Kundu}
 	\author { Dipendu Maity}
 	\affil{Department of Science and Mathematics,\\
 		Indian Institute of Information Technology Guwahati, Bongora, Assam-781\,015, India.\linebreak
 		\{arnab.kundu, dipendu\}@iiitg.ac.in}

\date{\today}

\maketitle

\begin{abstract}
If a map has $k$ transitivity classes of vertices that are subject to the action of the automorphism group, it is said to be $k$-uniform. The classification of $1$-uniform maps on the torus is known. In this article, we classify $2$-uniform maps on the torus up to isomorphism. Explicit formulas for the number of combinatorial types of these maps on number of vertices is obtained in terms of arithmetic functions in number theory, such as the divisor function. The asymptotic behaviour of these functions as number of vertices tends to infinity is also discussed and we obtained continuous functions which asymptotically served as upper and lower bounds. 
\end{abstract}

\noindent {\small {\em MSC 2010\,:} 52C20, 52B70, 51M20, 57M60.

 }

\section{Introduction}\label{intro}
Throughout the last few decades there have been many research going on with existence and classification of maps on surfaces. A surface is a connected 2-dimentional manifold.
A map intuitively means a Piecewise linear surface made by simplicial complexes or polygons. 
More precisely, a \textit{map} is a connected $2$-dimensional cell complex on a surface. Equivalently, it is a cellular embedding of a connected simple graph on a surface. 
For a map $\mathcal{K}$, let $V(\mathcal{K})$ be the vertex set of $\mathcal{K}$ and $u\in V(\mathcal{K})$. The faces containing  $u $ form a cycle (called the {\em face-cycle} at  $u $)  $C_u $ in the dual graph  of  $\mathcal{K} $. That is,  $C_u $ is of the form  $(F_{1,1}\mbox{-}\cdots\mbox{-}F_{1,n_1})\mbox{-}\cdots\mbox{-}(F_{k,1}\mbox{-}\cdots \mbox{-}F_{k,n_k})\mbox{-}F_{1,1}$, where  $F_{i,\ell} $ is a  regular $p_i $-gon for  $1\leq \ell \leq n_i $,  $1\leq i \leq k $,  $p_r\neq p_{r+1} $ for  $1\leq r\leq k-1 $ and  {$p_k\neq p_1 $}. In this case, the vertex $u$ is said to be of type $ [p_1^{n_1}, \dots, p_k^{n_k}]$ (addition in the suffix is modulo  $k $).
A map  $\mathcal{K} $ is said to be {\em $2$-semiequivelar} if  it contains only 2 types of vertices. If all vertices of a map is of same type then we call it \textit{semiequivelar}.
Two maps $\mathcal{K}_{1}$ and $\mathcal{K}_{2}$ are isomorphic if there exists a function $f ~:~ \mathcal{K}_{1}\rightarrow \mathcal{K}_{2}$ such that $f\mid _{V(\mathcal{K}_{1})} : V(\mathcal{K}_{1}) \rightarrow V(\mathcal{K}_{2})$ is a bijection and $f(\sigma)$ is a cell in $\mathcal{K}_{2}$ if and only if $\sigma$ is a cell in $\mathcal{K}_{1}$. In particular, if $\mathcal{K}_1 = \mathcal{K}_2$, then $f$ is called an $automorphism$. The \emph{automorphism group $\Aut(\mathcal{K})$} of $\mathcal{K}$ is the group consisting of automorphisms of $\mathcal{K}$.  A \emph{vertex-transitive} or \emph{1-uniform} map (or tiling) is a map (or tiling) on a closed surface (or on the plane) such that the automorphism group acts transitively on the set of vertices. 1-uniform tilings of $\mathbb{R}^2$ are called \textit{Archimedean tiling}. The $2$-uniform tilings (or map) are the generalization of vertex-transitive maps. A {\em $2$-uniform tiling} (or $map$) is a tiling (or a map) of the surface  having $2$ transitivity classes of vertices under the action of the automorphism group.  

Over the last few years there have been many results about maps and semi-equivelar maps that are highly symmetric. 
In particular, there has been recent interest in the study of maps using combinatorial, geometric, and algebraic approaches, with the topic of symmetries of maps receiving a lot of interest. We use topological and algebraic approach. To classify maps we use the concept of covering. 
\begin{definition}
Let $X$ and $Y$ be maps on surfaces $S$ and $S'$ respectively. A surjective function $\eta \colon X \to Y$ is called a $covering$ if it preserves adjacency and sends vertices, edges, faces of $X$ to vertices, edges, faces of $Y$ respectively.
\end{definition}
If $X$ and $Y$ has same type then for a  covering there exists a constant $n$ such that $n$ number of faces of $X$ maps to any face of $Y$. In this case we say the covering to be a \emph{n-sheeted covering}. Clearly $n$ is also the number of edges and vertices of $X$ maps to any edge or vertex of $Y$.

For a surface $S$ there exists a simply connected cover of $S$ called the \textit{universal cover}. This cover will be the hyperbolic plane($\mathbb{H}^2$), Euclidean plane($\mathbb{E}^2$) and the sphere($\mathbb{S}^2$) whenever $\chi(S)$(Euler characteristics of $S$) $<0, =0$ and $>0$ respectively. The following proposition gives an essential characterization of maps on surfaces.

\begin{proposition}
Let $X$ be a map on a surface $S$. Then $X$ can be obtained as an orbit space of a map on its universal cover $\widetilde{S}$ by action of some fixed point free discrete subgroup of  $\Aut(\widetilde{S})$.
\end{proposition}
This proposition tells that there is an one to one correspondence between maps on a surface $S$ with conjugacy classes of $\pi_1(S)$(the fundamental group of $S$) in $\Aut(\widetilde{S})$. 
In context of coverings there is an interest to classify covers of certain maps up to isomorphism. In particular,  there is also much interest in finding minimal regular covers of different families of maps and polytopes (see \cite{HW2012, MPW2013, pw2011}). In \cite{drach:2015}, Drach et al. constructed the minimal rotary cover(a cover which is either regular or chiral) of any equivelar toroidal map. Then, they have extended their idea to toroidal maps that are no longer equivelar, and constructed minimal toroidal covers of the Archimedean toroidal maps with maximal symmetry (see in \cite{drach:2019}), called these covers almost regular; they will no longer be regular (or chiral), but instead will have the same number of flag orbits as their associated tessellation in Euclidean plane. In \cite{KM1,KM2,KM3} authors proved existence of more symmetric covers of Archimedean toroidal map from three points of view; vertices, edges and flags. They also enumerated $n$ sheeted covers of a given map.

Maps on compact surfaces are orbit space of the universal cover of that surface by a discrete subgroups of the symmetry group. Two maps are isomorphic if and only if corresponding groups are conjugate in the symmetry group of the universal cover.
In this context Mednykh gave a formula in \cite{Med} for the number of conjugacy
classes of subgroups of a given index in a finitely generated group. As application of
this result a simple proof of the formula for the number of non-equivalent coverings
over surface (orientable or not, bordered or not) is given.
In \cite{LM1990, LM1991} Lu\v{c}i\'{c} and Moln\'{a}r discussed about the geometrical structure of fundamental domain and gave an algorithm for complete enumeration of uniform tilings of any complete, simply connected, two dimensional Riemmanian manifold of constant curvature in particular of all non congruent uniform tilings of hyperbolic plane. This method based on the classification of planer discontinuous groups which also enumerates all combinatorially different fundamental domains for any given planer discontinuous group.


As of now, the only significant study on the classification of maps on the sphere, projective plane, torus, Klein bottle, non-orientable surfaces of genus $3$ and orientable surfaces of genus $2, 3, 4$ has been done by many researchers. More precisely, the well known eleven types of semiregular tiling of the plane are the examples of vertex transitive semiequivelar maps on Euclidean plane. All semiequivelar and vertex-transitive maps on the 2-sphere are known. These are the boundaries of Platonic and Archimedean solids and two infinite families of types (namely, of types $[4^2, n^1]$ and $[3^3,m^1]$ for $4\neq n\geq 3$, $m\geq 4$) \cite{DM2022, GS1977}. 
Similarly, there are infinitely many types of semiequivelar and vertex-transitive on the real projective plane \cite{Ba1991, DM2022}. Torus is a quotient space of Euclidean plane. One can construct infinitely many semiequivelar maps on this closed surface as a quotient of one of the eleven semiregular tilings by a group $\Gamma$ generated by two translations \cite{Su2011t, Su2011kb}. Karabas and Nedela have obtained the census of Archimedean maps on orientable surface of Euler characteristic $-2$. In 2012, they have classified all the Archimedean maps on orientable surface of Euler characteristic $-2$, $-4$, $-6$ \cite{kn2012}. Kurth, Negami, Brehm and K\"{u}hnel  have studied all equivelar maps of type $[3^{6}]$, $[4^{4}]$, $[6^{3}]$ on the torus and given complete classification independently \cite{alt1973, brehm2008, kurth:1986, negami1983}. 

Here we fix our surface to be torus and address the following question: 

\begin{question}\label{ques}
Does there exist any map $2$-uniform map on the torus? If so, how many? 
\end{question}

In this context, we know the following results.


\begin{proposition} (\cite{alt1973, MU2018}) \label{prop-equi-maps} Let $X$ be an Archimedean toroidal map. Then we have classification of maps of type $X$ on a given number of vertices up to isomorphism. 
\end{proposition}

\begin{proposition} (\cite{M2022}) Let $X$ be an Archimedean toroidal map.  Then we have classification of vertex-transitive maps ($1$-uniform maps) of type $X$ on a given number of vertices up to isomorphism. 
\end{proposition}

Thus, we know the classification of $1$-uniform maps on the torus. In this article, we classify $2$-uniform toroidal maps. In particular, we prove the followings:



\begin{theorem}\label{thm-1}
{\rm (a)} Let $X$ be a 2-uniform toroidal map of type $\ell$ with $v$ vertices. We give explicit functions ($\Phi_{\ell})$ on number of vertices ($v$) to enumerate such possible maps. In Table \ref{tabl1} we give $\Phi_{\ell}(v)$ for some $v$.\\
{\rm (b)} We devise a way to obtain the list of 2-uniform toroidal maps (up to isomorphism). \\
{\rm (c)} We also discuss some continuous functions which are asymptotically bounds $\Phi_{\ell}$ for each $\ell$.
\end{theorem}

This paper organised as follows.
In Section \ref{sec:examples} we give the pictures of the 2-uniform and Archimedean tilings of the plane and the Table \ref{tabl1} where we give explicit number of $2$-uniform toroidal maps for given number of vertex. These values are calculated using SageMath. 
In Section \ref{background} we give further background on association of toroidal maps with non singular matrices and will talk about their uniqueness. We also give a brief description of our approach.
Then in Section \ref{proof} we enumerate the 2-uniform toroidal maps on their number of vertices and give explicit functions for the number of combinatorial types. Finally we discuss about asymptotic behaviour some of the obtained functions in Section \ref{ab} and remarked that the number of vertices $v$ is given such that the number of 2-uniform maps on torus up to isomorphism with $v$ vertices is larger than $v$ itself.

\section{Examples}\label{sec:examples}
Suppose $E$ is a tiling on $\mathbb{R}^2$ and $X = E/K$ be a toroidal map where $K\le \Aut(E)$. Then, the number of vertex orbits of $X$ under the action of automorphism group of $X$ will be at least the number of vertex orbits of $E$ under the automorphism group of $E$ since any automorphism of $X$ can be lifted to an automorphism of $E$. Hence, if $X$ has two vertex orbits then number of vertex orbits of $E$ will be at most two. Thus, $E$ must belong to either $2$-uniform tilings of the plane or the $1$-uniform tilings of the plane. 
We know from \cite{GS1977, GS1981, Otto1977} that there are {twenty} 2-uniform tilings and 7 Archimedean tilings which are not always vertex transitive on torus.
Here are these tilings. The vertex types of them are written in bracket.

\begin{figure}[H]
    \centering
    \subfloat[{$E_1 ~ ([3^{6};3^4,6^1])$}]{
            \includegraphics[height=3.72cm, width= 3.70cm]{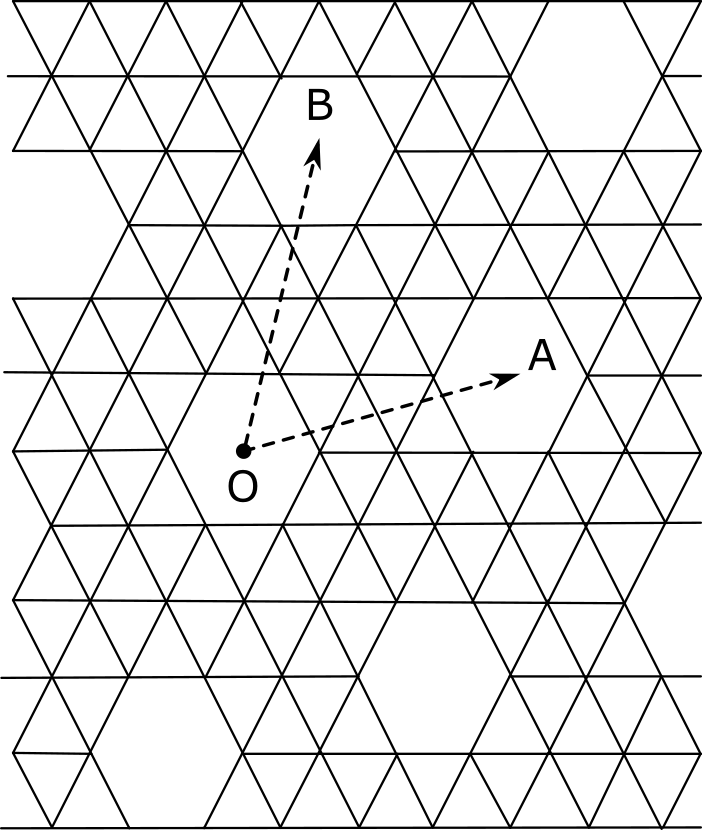}\label{E_1}
    }
    \subfloat[{$E_2$~($[3^{6};3^4,6^1]$)}]{
        \includegraphics[height=3.72cm, width= 3.70cm]{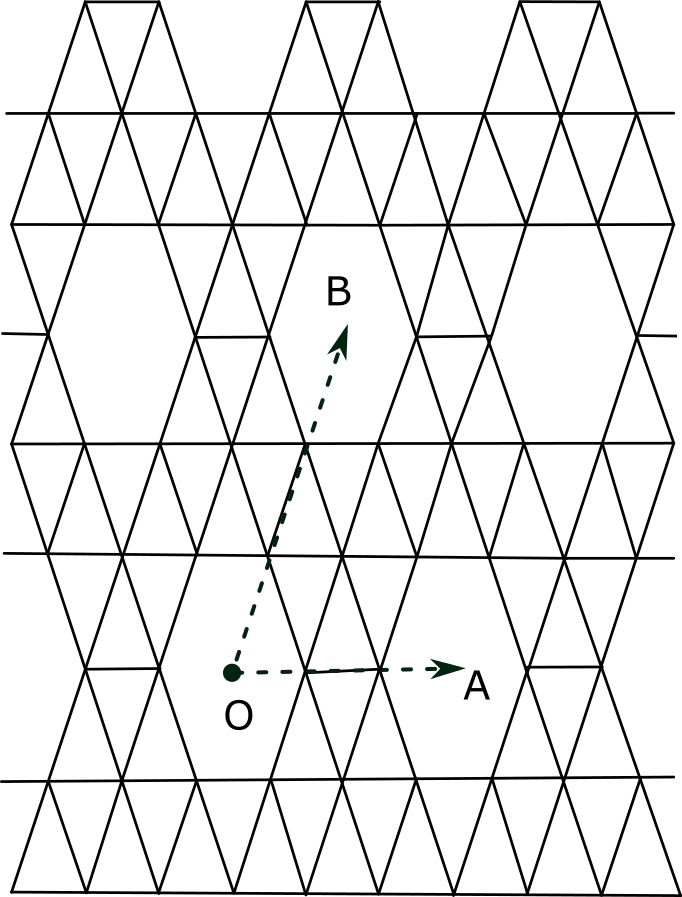}\label{E_2}
    }     
    \subfloat[{$E_3 ~ ([3^6;3^3,4^{2}])$}]{
            \includegraphics[height=3.72cm, width= 3.72cm]{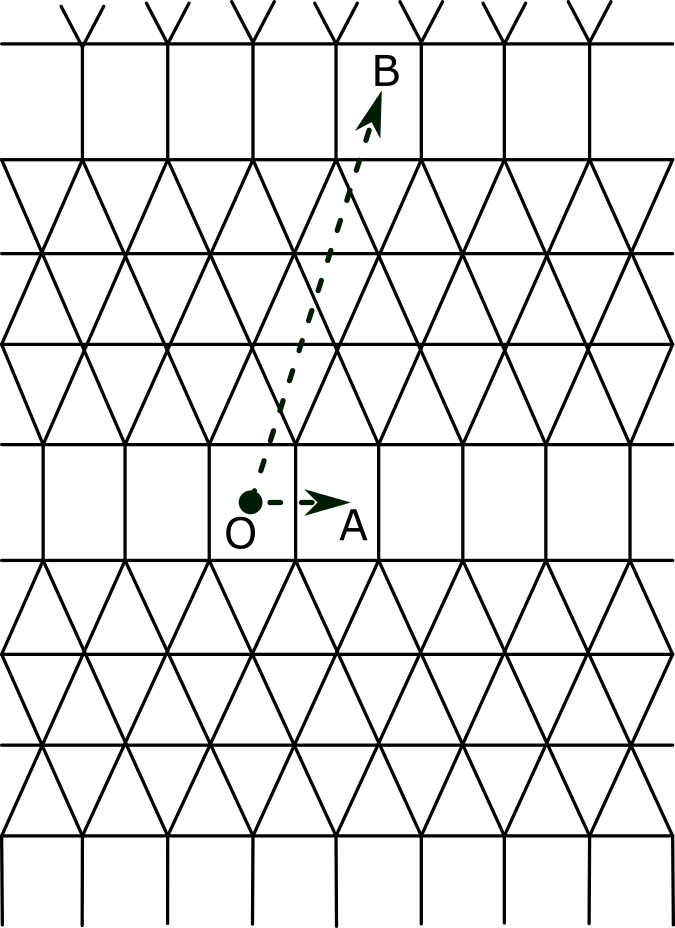}\label{E_3}
            
    }
    \subfloat[{$E_4$~($[3^{6};3^3,4^2]$)}]{
        \includegraphics[height=3.72cm, width= 3.72cm]{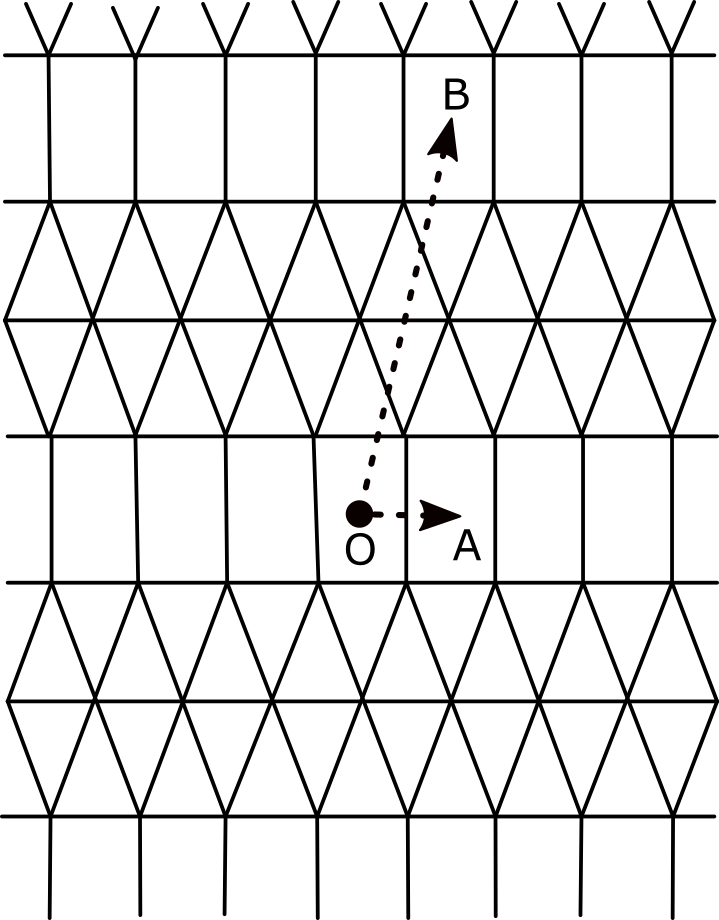}\label{E_4}
    }\\ 
    
      
     \subfloat[{$E_5$~($[3^{6};3^2,4^1,3^1,4^1]$)}]{
        \includegraphics[height=3.72cm, width= 3.8cm]{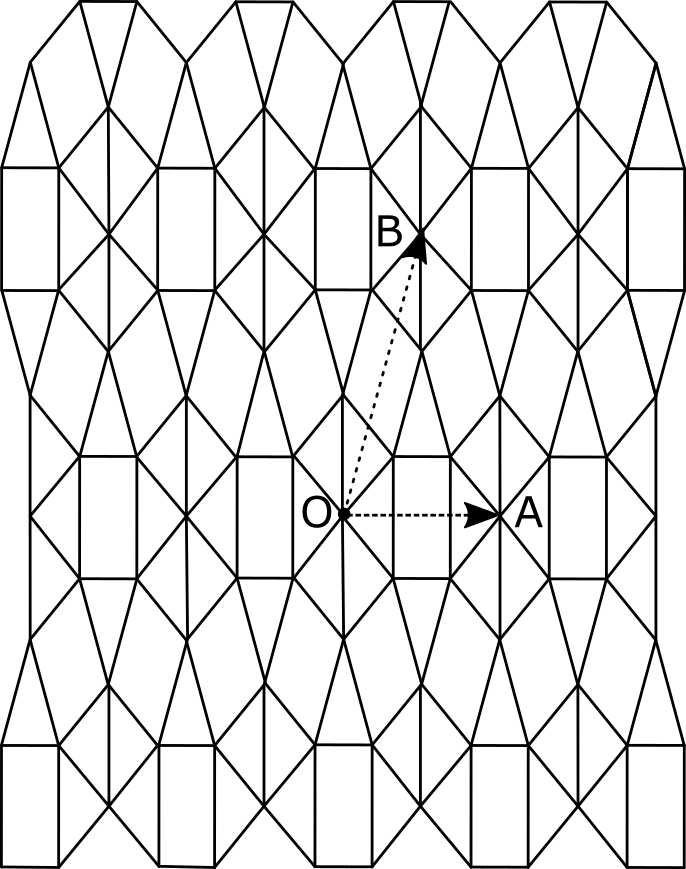}\label{E_5}
    }   
    \subfloat[{$E_6$~($[3^{6};3^2,4^1,12^1]$)}]{
        \includegraphics[height=3.72cm, width= 3.64cm]{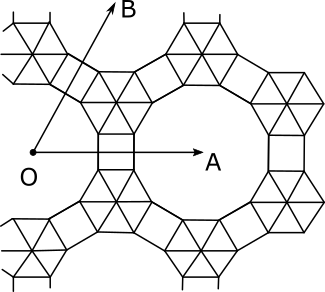}\label{E_6}
    }  
     \subfloat[{$E_7$~($[3^{6};3^2,6^2]$)}]{
        \includegraphics[height=3.72cm, width= 3.70cm]{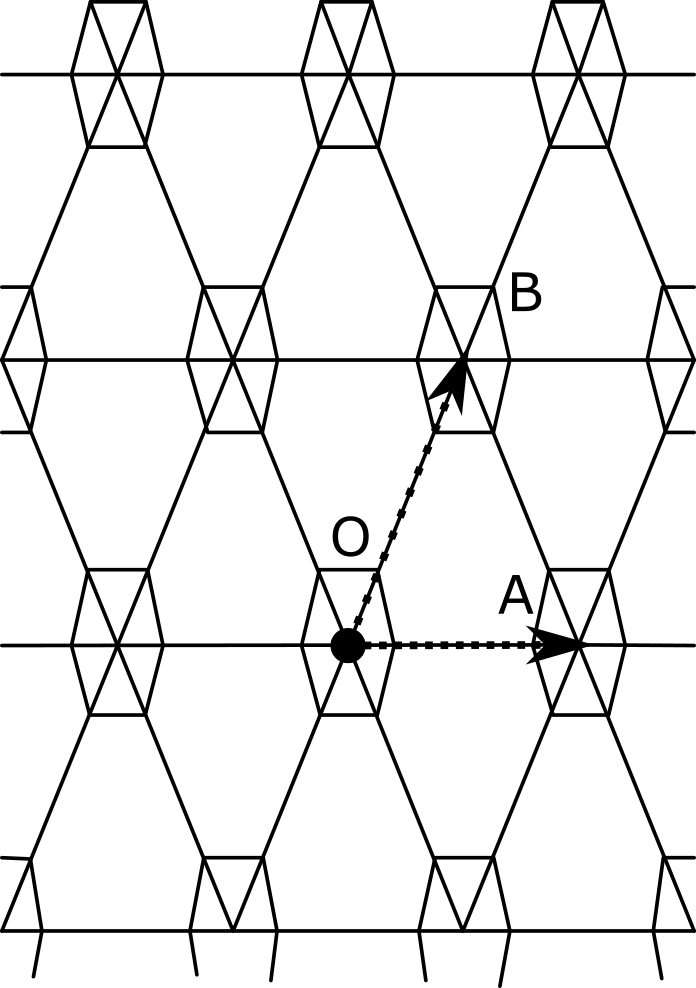}\label{E_7}
    } 
     \subfloat[{$E_8$~($[3^{2},6^2;3^4,6^1]$)}]{
        \includegraphics[height=3.72cm, width= 3.72cm]{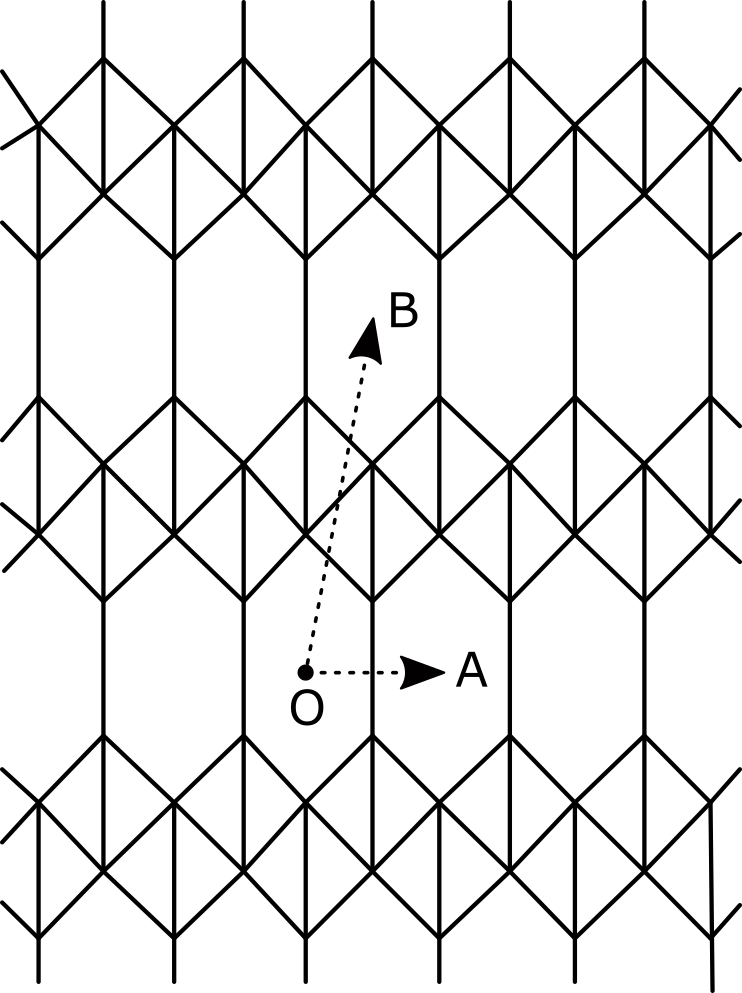}\label{E_8}
    } 
    \\ 
     \end{figure}
    
      \begin{figure}[H]
      \centering
    
     \subfloat[{$E_9$~($[3^{3},4^2;3^2,4^1,3^1,4^1]$)}]{
        \includegraphics[height=3.2cm, width= 3.72cm]{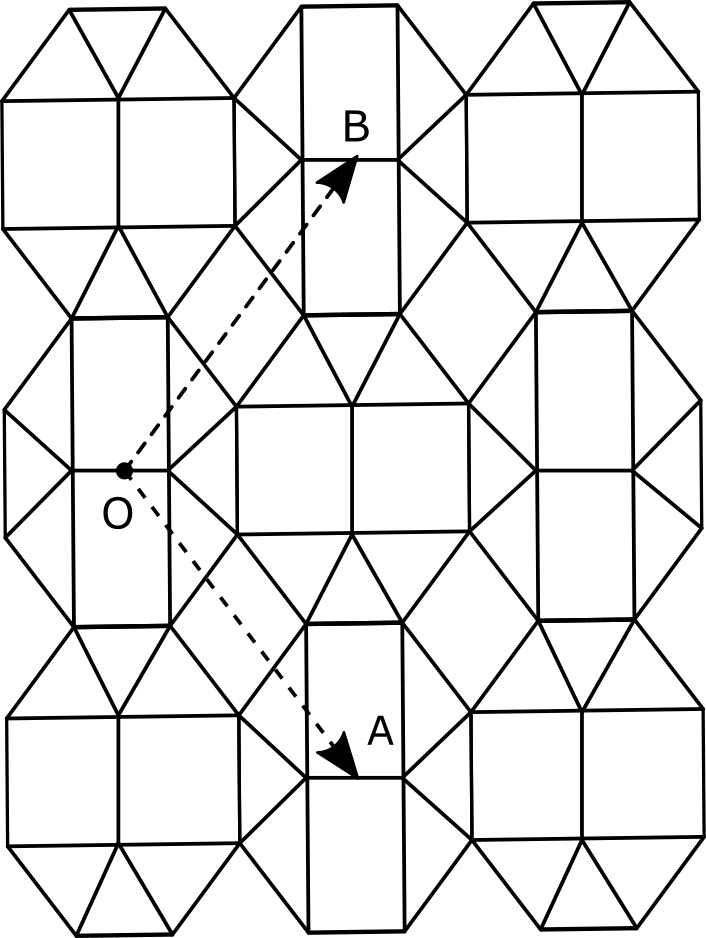}\label{E_9}
    } 
     \subfloat[{$E_{10}$~($[3^{3},4^2;3^2,4^1,3^1,4^1]$)}]{
        \includegraphics[height=3.2cm, width= 3.70cm]{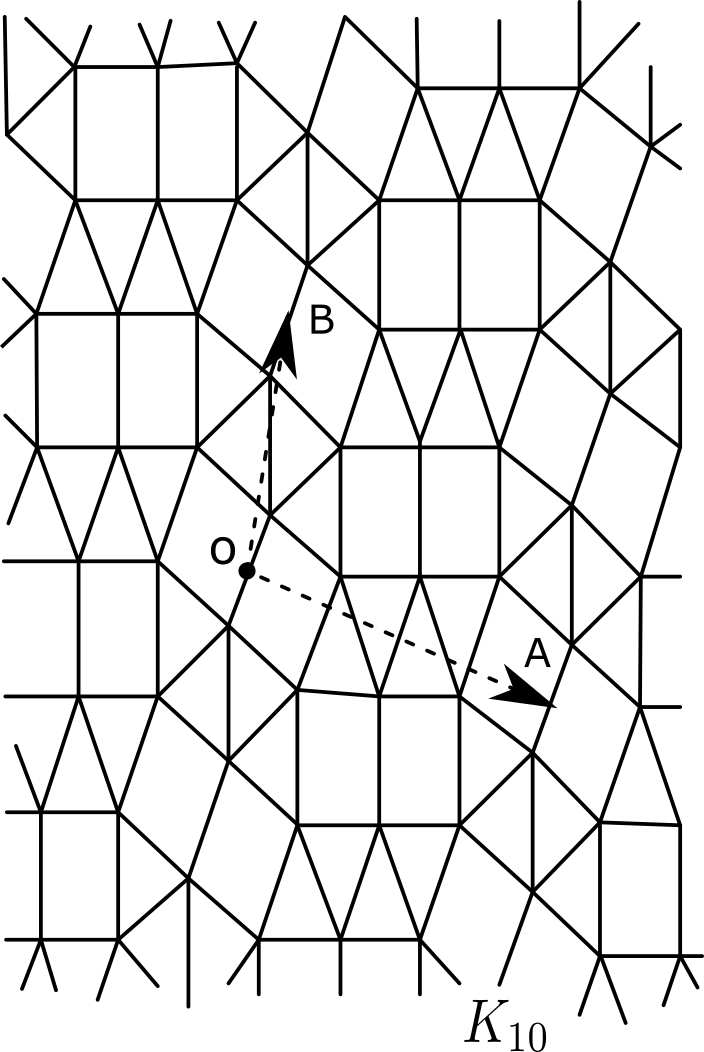}\label{E_10}
    } 
     \subfloat[{$E_{11}$~($[3^{3},4^2;3^1,4^1,6^1,4^1]$)}]{
        \includegraphics[height=3.2cm, width= 3.70cm]{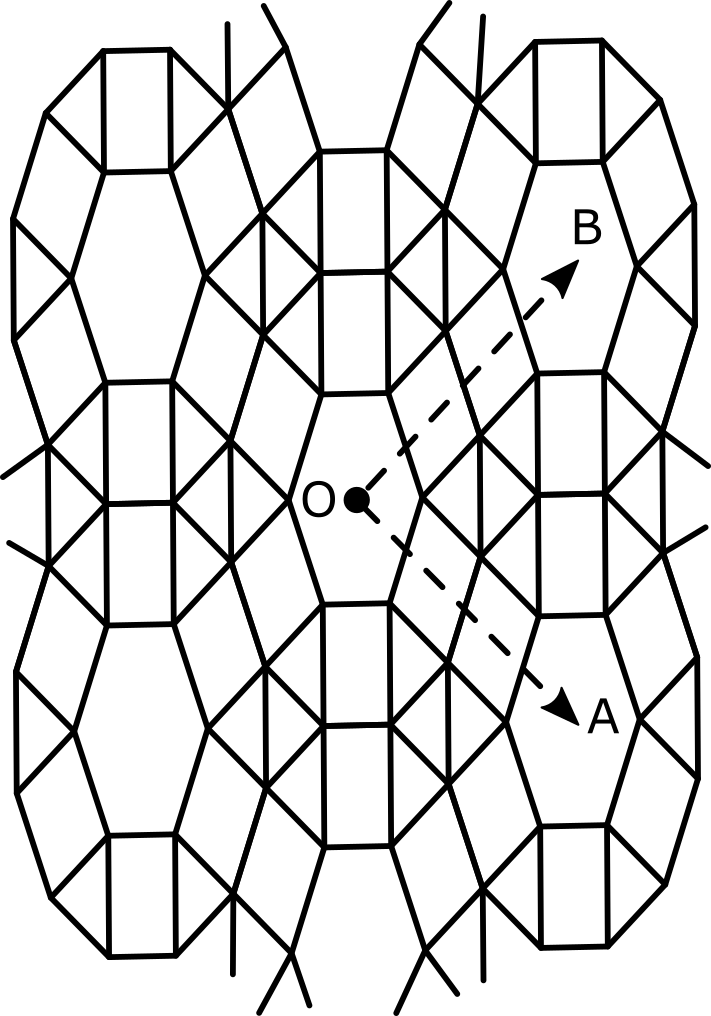}\label{E_11}
    } 
     \subfloat[{$E_{12}$~($[4^{4};3^3,4^2]$)}]{
        \includegraphics[height=3.2cm, width= 3.72cm]{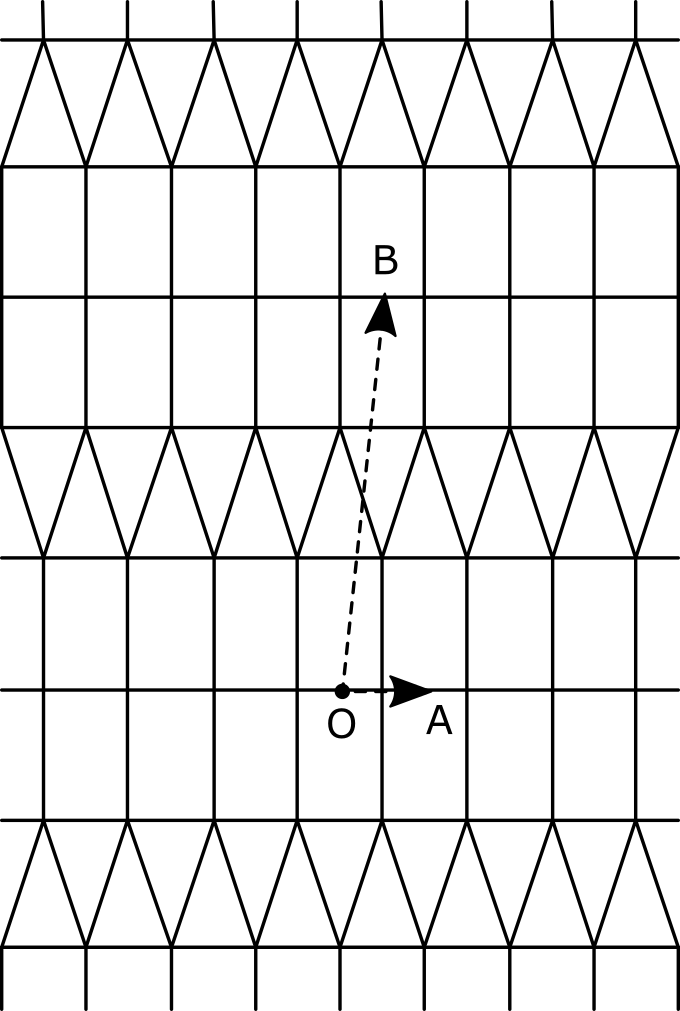}\label{E_12}
    } \\ 
    
    
    \subfloat[{$E_{13}$~($[4^4,4^2;3^3,4^2]$)}]{
        \includegraphics[height=3.2cm, width= 3.72cm]{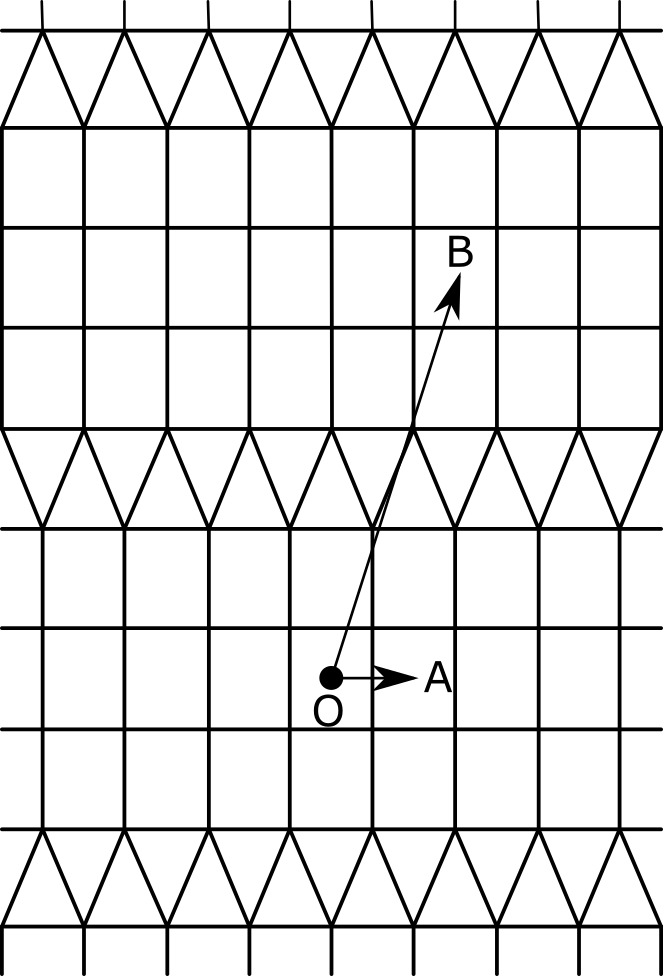}\label{E_13}
    }
     \subfloat[{$E_{14}$~($[3^1,4^1,6^1,4^1;3^2,4^1,3^1,4^1]$)}]{
        \includegraphics[height=3.2cm, width= 3.72cm]{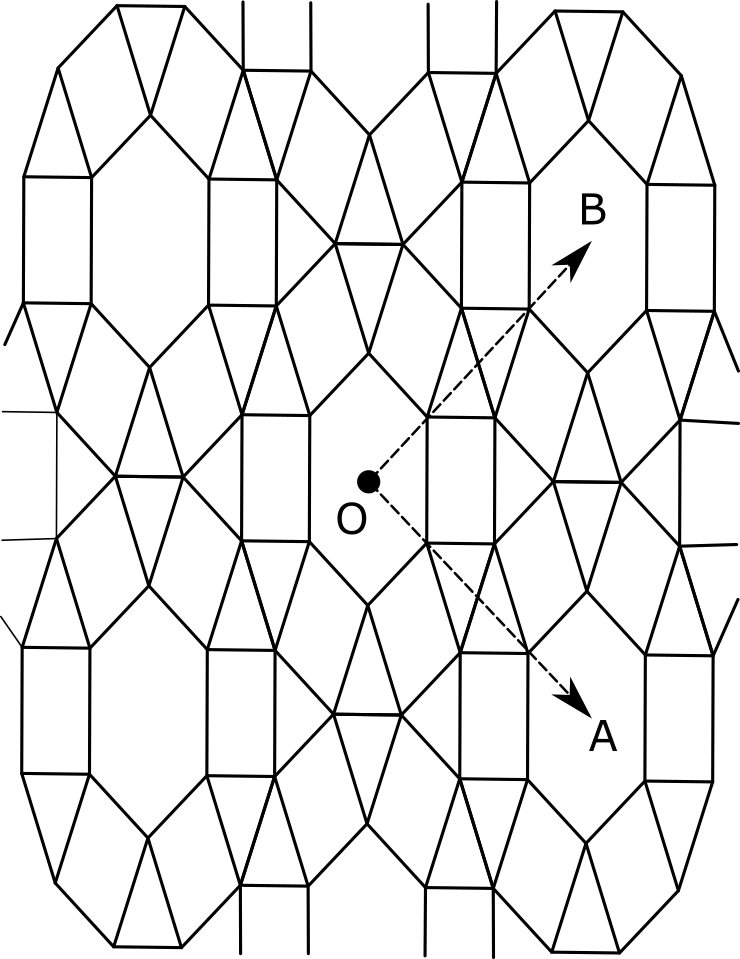}\label{E_14}
    }
    \subfloat[{$E_{15}$~($[3^2,6^2;3^1,6^1,3^1,6^1]$)}]{
        \includegraphics[height=3.2cm, width= 3.70cm]{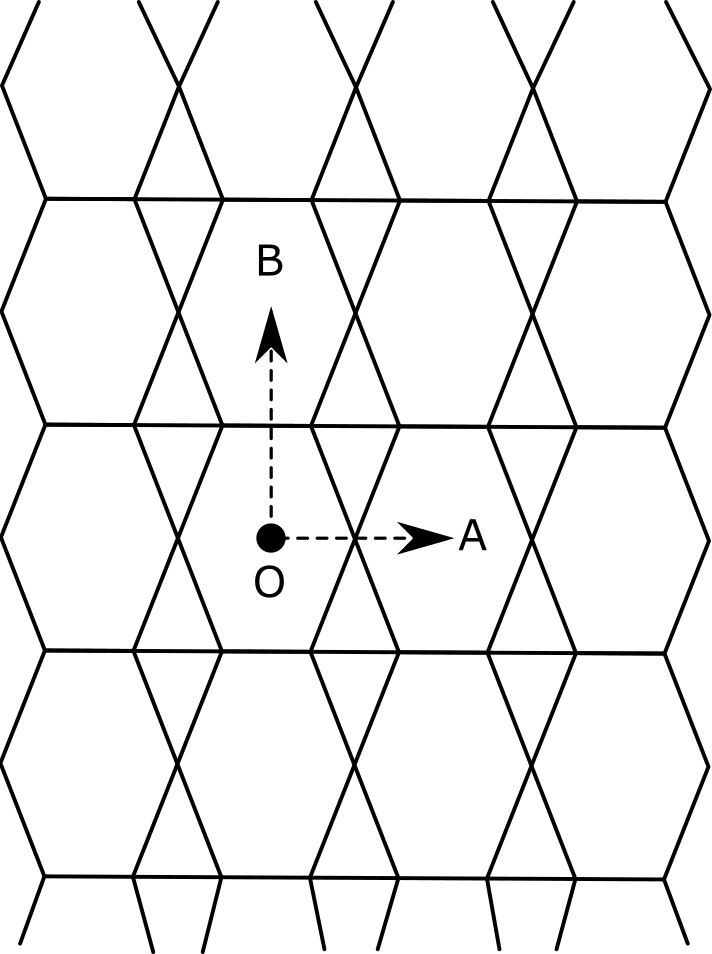}\label{E_15}
    }
    \subfloat[{$E_{16}$~($[12^2,3^1;3^1,4^1,3^1,12^1]$)}]{
        \includegraphics[height=3.2cm, width= 3.70cm]{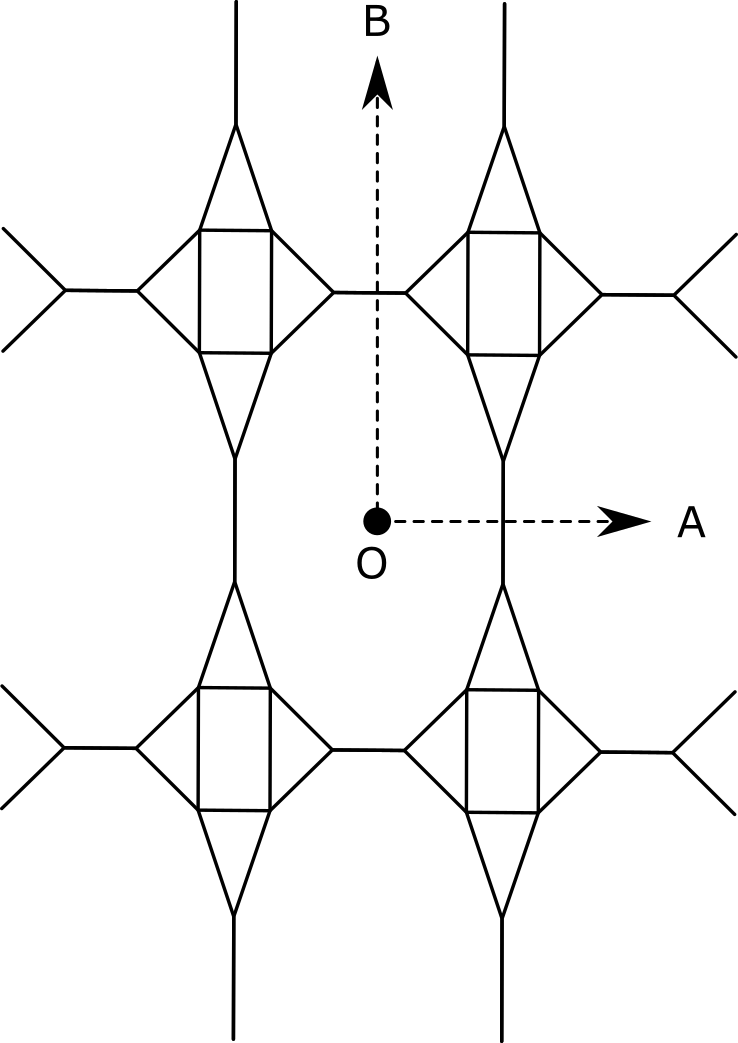}\label{E_16}
    }\\
    
    

    
    \subfloat[{$E_{17}$~($[3^1,4^1,6^1,4^1;3^1,4^2,6^1]$)}]{
        \includegraphics[height=3.2cm, width= 3.72cm]{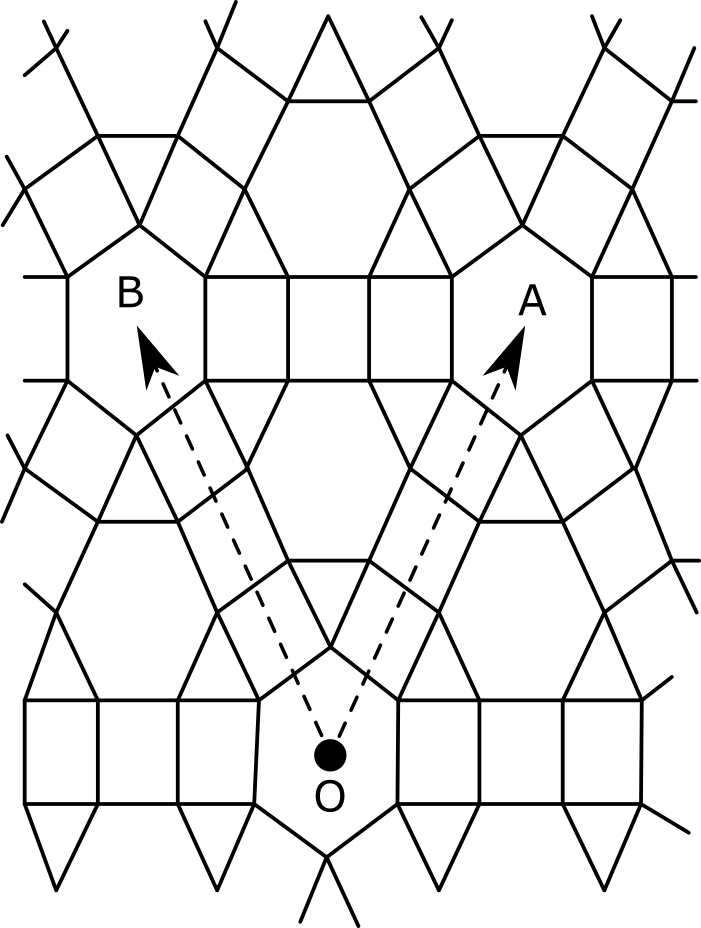}\label{E_17}
    }
    \subfloat[{$E_{18}$~($[3^1,4^2,6^1;3^1,6^1,3^1,6^1]$)}]{
        \includegraphics[height=3.2cm, width= 3.72cm]{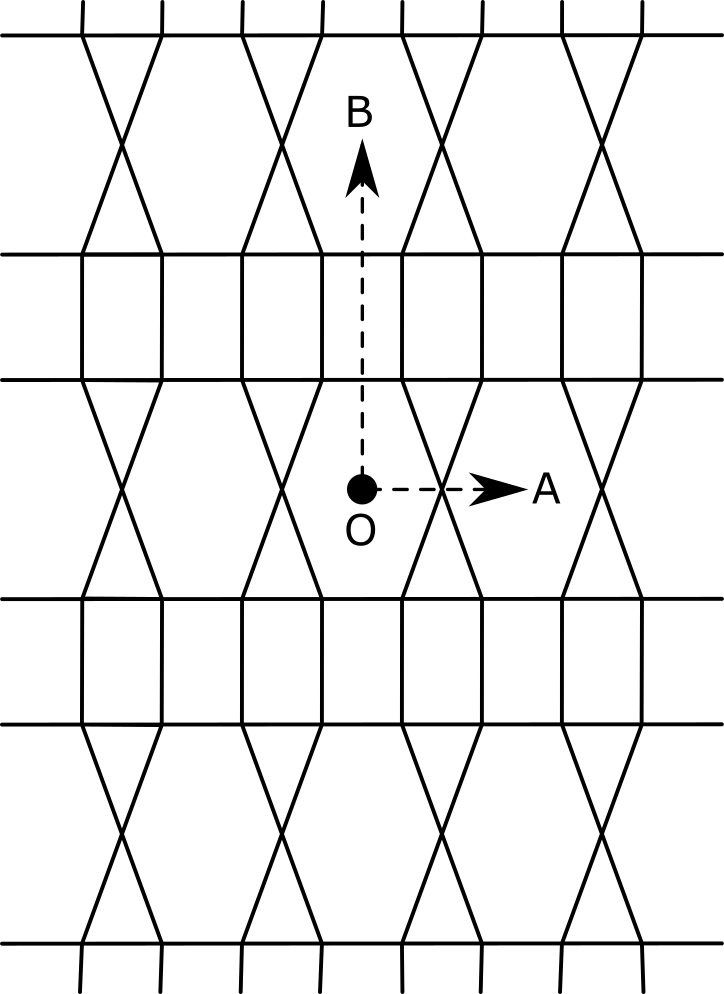}\label{E_18}
    }
    \subfloat[{$E_{19}$~($[3^1,6^1,3^1,6^1;3^1,4^2,6^1]$)}]{
        \includegraphics[height=3.2cm, width= 3.72cm]{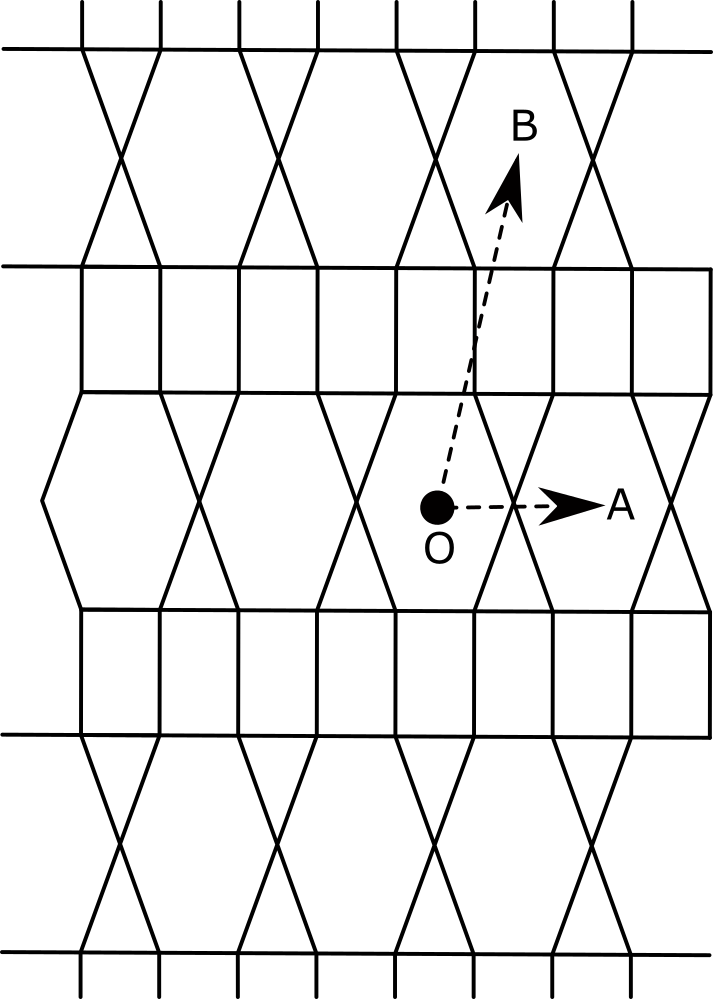}\label{E_19}
    }
    \subfloat[{$E_{20}$~($[4^1,6^1,12^1;3^1,4^1,6^1,4^1]$)}]{
        \includegraphics[height=3.2cm, width= 3.72cm]{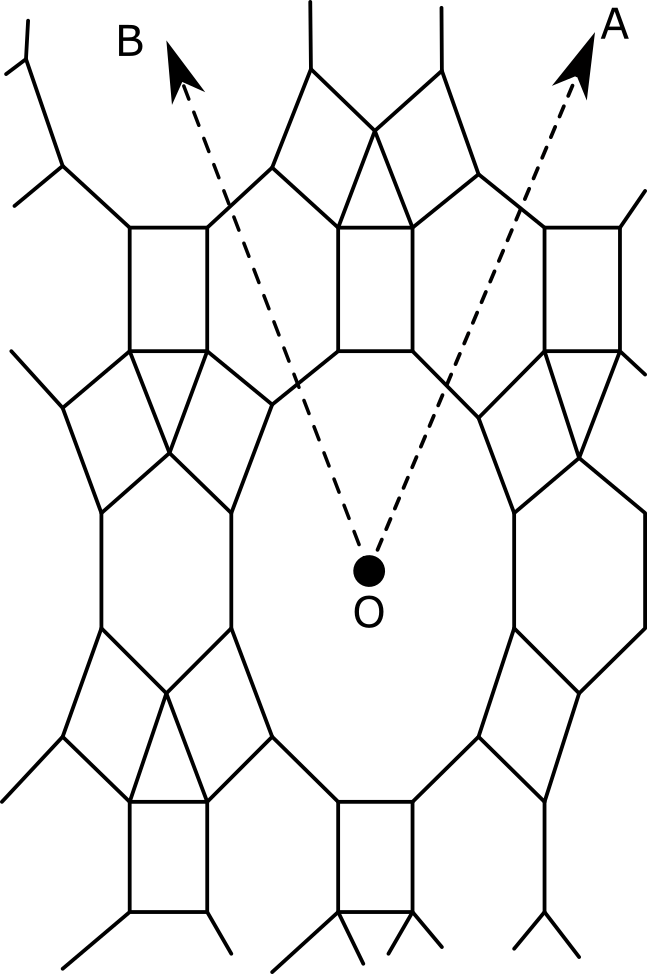}\label{E_20}
    }
    \caption{2-uniform tilings of the plane.}\label{Fig-1}
    \end{figure}
    
   \vspace{-1cm}
    
    \begin{figure}[H]
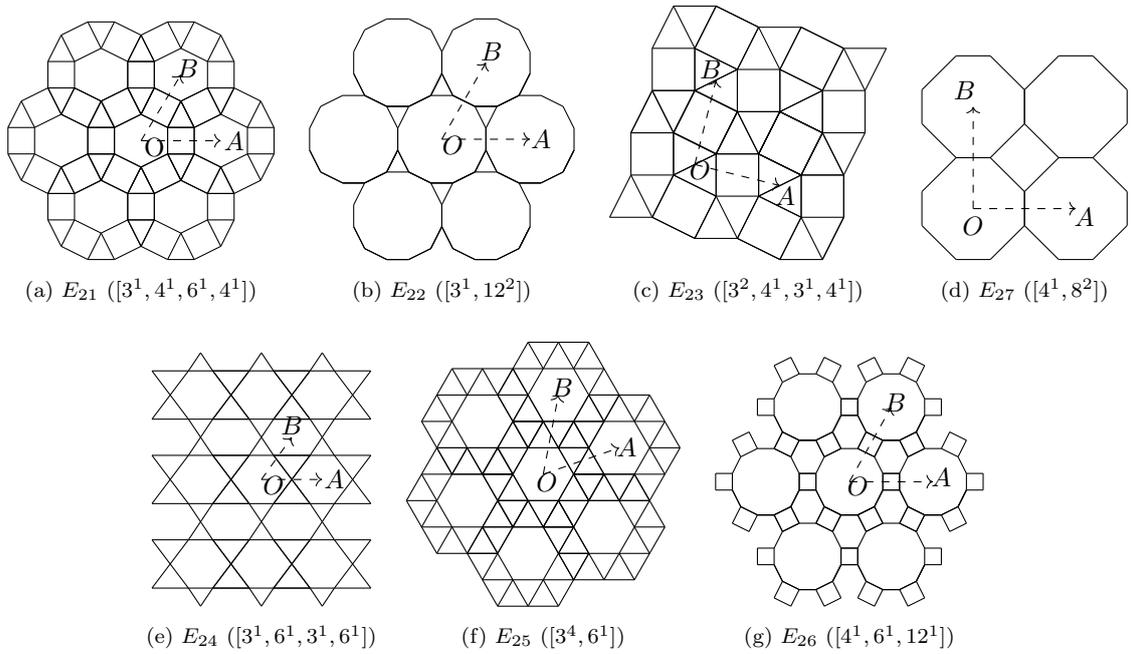
\label{fig-E_9} 
    \centering
    \subfloat[{$E_{21}$~($[3^1,4^1,6^1,4^1]$)}]{


    }
    
     \caption{Archimidean tilings of the plane.}\label{fig-2}
    \end{figure}
     
\begin{landscape}
\begin{table}
\centering
{\small

\begin{tabular}{||c|c|| c |c||c |c||c |c||c |c||c |c||c |c||c |c||c |c||c |c||c|c ||c|c ||c|c ||}
    \hline
        v  &  $\Phi_{1}(v)$ &v  &  $\Phi_{2}(v)$ &v  &  $\Phi_{5}(v)$ &v  &  $\Phi_{6}(v)$ &v  &  $\Phi_{7}(v)$ &v  &  $\Phi_{8}(v)$ &v  &  $\Phi_{9}(v)$ &v  &  $\Phi_{10}(v)$&v  &  $\Phi_{11}(v)$  \\ \hline
        12 & 1 & 8 & 1 & 7 & 1 & 14 & 1 & 7 & 1 & 4 & 1 & 12 & 1 & 8 & 1 & 12 & 1\\ \hline
        24 & 0 & 16 & 0 & 14 & 0 & 28 & 0 & 14 & 0 & 8 & 2 & 24 & 3 & 16 & 3 & 24 & 0\\ \hline
        36 & 1 & 24 & 1 & 21 & 1 & 42 & 1 & 21 & 1 & 12 & 3 & 36 & 2 & 24 & 2 & 36 & 1\\ \hline
        48 & 1 & 32 & 1 & 28 & 1 & 56 & 1 & 28 & 1 & 16 & 5 & 48 & 5 & 32 & 5 & 48 & 1\\ \hline
        60 & 0 & 40 & 0 & 35 & 0 & 70 & 0 & 35 & 0 & 20 & 4 & 60 & 2 & 40 & 2 & 60 & 0\\ \hline
        72 & 0 & 48 & 0 & 42 & 0 & 84 & 0 & 42 & 0 & 24 & 7 & 72 & 6 & 48 & 6 & 72 & 0\\ \hline
        84 & 2 & 56 & 1 & 49 & 1 & 98 & 0 & 49 & 1 & 28 & 5 & 84 & 2 & 56 & 2 & 84 & 1\\ \hline
        96 & 0 & 64 & 0 & 56 & 0 & 112 & 0 & 56 & 0 & 32 & 10 & 96 & 7 & 64 & 7 & 96 & 0\\ \hline
        108 & 1 & 72 & 1 & 63 & 1 & 126 & 1 & 63 & 1 & 36 & 8 & 108 & 3 & 72 & 3 & 108 & 1\\ \hline
        120 & 0 & 80 & 0 & 70 & 0 & 140 & 0 & 70 & 0 & 40 & 10 & 120 & 6 & 80 & 6 & 120 & 0\\ \hline
    \end{tabular}

\vspace{1cm}

\begin{tabular}{||c|c|| c |c||c |c||c |c||c |c||c |c||c |c||c |c||c |c||c |c||c|c ||c|c ||c|c ||}
    \hline
        v  &  $\Phi_{12}(v)$ &v  &  $\Phi_{13}(v)$ &v  &  $\Phi_{14}(v)$& v  &  $\Phi_{15}(v)$&v  &  $\Phi_{16}(v)$ &v  &  $\Phi_{17}(v)$ &v  &  $\Phi_{18}(v)$ &v  &  $\Phi_{19}(v)$ &v  &  $\Phi_{20}(v)$ \\ \hline
      3 & 1 & 4 & 1 & 12 & 1 & 3 & 1 & 8 & 1 & 18 & 1 & 5 & 1 & 5 & 1 & 18 & 1\\ \hline
      6 & 2 & 8 & 2 & 24 & 0 & 6 & 3 & 16 & 1 & 36 & 0 & 10 & 3 & 10 & 1 & 36 & 0\\ \hline
      9 & 3 & 12 & 3 & 36 & 1 & 9 & 3 & 24 & 0 & 54 & 1 & 15 & 2 & 15 & 2 & 54 & 1\\ \hline
      12 & 5 & 16 & 5 & 48 & 1 & 12 & 6 & 32 & 1 & 72 & 1 & 20 & 5 & 20 & 3 & 72 & 1\\ \hline
      15 & 4 & 20 & 4 & 60 & 0 & 15 & 4 & 40 & 1 & 90 & 0 & 25 & 2 & 25 & 2 & 90 & 0\\ \hline
      18 & 7 & 24 & 7 & 72 & 0 & 18 & 9 & 48 & 0 & 108 & 0 & 30 & 6 & 30 & 2 & 108 & 0\\ \hline
      21 & 5 & 28 & 5 & 84 & 1 & 21 & 5 & 56 & 0 & 126 & 0 & 35 & 2 & 35 & 2 & 126 & 0\\ \hline
      24 & 10 & 32 & 10 & 96 & 0 & 24 & 11 & 64 & 1 & 144 & 0 & 40 & 7 & 40 & 5 & 144 & 0\\ \hline
      27 & 8 & 36 & 8 & 108 & 1 & 27 & 8 & 72 & 1 & 162 & 1 & 45 & 3 & 45 & 3 & 162 & 1\\ \hline
      30 & 10 & 40 & 10 & 120 & 0 & 30 & 12 & 80 & 1 & 180 & 0 & 50 & 6 & 50 & 2 & 180 & 0\\ \hline
    \end{tabular}
    
    \vspace{1cm}
    
\begin{tabular}{||c|c|| c |c||c |c||c |c||c |c||c |c||c |c||c |c||c |c||c |c||c|c ||c|c ||c|c ||}
    \hline
        v  &  $\Phi_{21}(v)$ &v  &  $\Phi_{22}(v)$ &v  &  $\Phi_{23}(v)$ &v  &  $\Phi_{24}(v)$ &v  &  $\Phi_{26}(v)$ &v  &  $\Phi_{27}(v)$  \\ \hline
      6 & 0 & 6 & 0 & 4 & 0 & 3 & 0 & 12 & 0 & 4 & 0\\ \hline
      12 & 1 & 12 & 1 & 8 & 0 & 6 & 1 & 24 & 0 & 8 & 0\\ \hline
      18 & 1 & 18 & 1 & 12 & 1 & 9 & 1 & 36 & 0 & 12 & 0\\ \hline
      24 & 2 & 24 & 2 & 16 & 1 & 12 & 2 & 48 & 0 & 16 & 0\\ \hline
      30 & 2 & 30 & 2 & 20 & 1 & 15 & 2 & 60 & 0 & 20 & 0\\ \hline
      36 & 2 & 36 & 2 & 24 & 2 & 18 & 2 & 72 & 0 & 24 & 2\\ \hline
      42 & 2 & 42 & 2 & 28 & 2 & 21 & 2 & 84 & 1 & 28 & 2\\ \hline
      48 & 5 & 48 & 5 & 32 & 3 & 24 & 5 & 96 & 0 & 32 & 2\\ \hline
      54 & 2 & 54 & 2 & 36 & 3 & 27 & 2 & 108 & 0 & 36 & 4\\ \hline
      60 & 2 & 60 & 2 & 40 & 3 & 30 & 2 & 120 & 0 & 40 & 4\\ \hline
    \end{tabular}
  
}

\caption{Enumeration of $2$-uniform toroidal maps}
\label{tabl1}
\end{table}

\end{landscape}

\section{Preliminaries}\label{background}
In this section we provide definitions and results which are required to proof our main theorems.

Let $X$ be a toroidal map. Then $X$ is obtained as an orbit space of some tiling $E$ on $\mathbb{R}^2$. Then, we can write $X=\frac{E}{K}$ where $K\leq\Aut(E)$ be discrete and does not have any fixed points. So, $K$ contains only translations and glide reflections. Being a torus, $X$ is orientable, thus $K\leq H$ where $H$ be the group of translational symmetries of $E$. Let us choose some origin $O$ in $E$. Let $H=\langle \gamma,\,\delta\rangle$, where $\gamma: z\mapsto  z+A$ and $\delta: z\mapsto  z+B$ where $A$ and $B$ are two linearly independent translation vectors of $E$ originating from $O$. Then, $K=\langle w_1, w_2\rangle$ where $w_1=\gamma^a\circ \delta^b$ and $w_2=\gamma^c\circ \delta^d$ for some $a,\,b,\,c,\,d\in\mathbb{Z}$. Clearly, $w_1:z\mapsto z+(aA+bB)$ and $w_2: z\mapsto z+(cA+dB)$.
Then $X$ is associated with $\begin{bmatrix} a & c\\ b & d\end{bmatrix}$. 

The matrix representation is unique for the fundamental domain of a map but not for the map itself, since a map can be represented by different fundamental domains. The transformation from one matrix to another corresponds to shearing the fundamental domain. The rotations and reflections will not change the combinatorial type of the map but the associated matrix will be multiplied by some matrix from the left. We note this in the following proposition.

\begin{proposition}\cite{brehm2008}
Let $X_1=E/K_1$ and $X_2=E/K_2$ be two toroidal maps of same type. Let $M_1$ and $M_2$ be their associated matrices. Then $X_1$ and $X_2$ will be isomorphic if and only if there exists matrices $A\in S_o$ and $B\in GL(2,\ZZ)$ such that $M_2 = AM_1B$. Here $S_o$ is the stabilizer of the origin by the action of $\Aut(E)$ on $E$.
\end{proposition}

We call two maps isomorphic if they are isomorphic as maps (see Section \ref{intro}).
Two maps are equal if the orbits of $\mathbb{R}^2$ under the action of the corresponding
groups are identical as sets. From \cite{kurth:1986} we have the following.

\begin{proposition}
Let $A$ and $B$ be $2 \times 2$ matrices with integral entries. The maps
corresponding to them are equal if and only if there exists a unimodular matrix $U$ with
integral entries satisfying $A U = B$.
\end{proposition}

The matrix representation becomes unique if the matrices of type $\begin{bmatrix}a&0\\b&d\end{bmatrix}$ with restriction $a,b,d>0$ and $0\le b<d$ \cite{kurth:1986}. This representation of integer matrices are called hermite normal form(HNF) \cite{wiki}.


For example, the maps red and green in Fig. \ref{eg} bellow has associated matrix $\begin{bmatrix} 4&0\\3&3\end{bmatrix}$ and $\begin{bmatrix} 4&0\\0&3\end{bmatrix}.$ These two maps are isomorphic to $M$ but not equal.

\begin{figure}[H]
    \centering
    \subfloat[]{
            \includegraphics[scale=0.48]{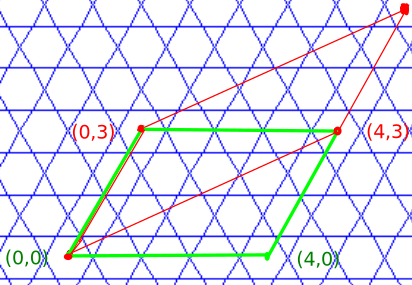}
    }
    \subfloat[$M$]{
        \begin{tikzpicture}[scale=0.7]
{\small
		\node  (0) at (-6, 4) {};
		\node  (1) at (-4, 4) {};
		\node  (2) at (-2, 4) {};
		\node  (3) at (-5, 4) {};
		\node  (4) at (-3, 4) {};
		\node  (5) at (-1, 4) {};
		\node  (6) at (0, 4) {};
		\node  (7) at (-6, 3) {};
		\node  (8) at (-6, 2) {};
		\node  (9) at (-5, 2) {};
		\node  (10) at (-4, 2) {};
		\node  (11) at (-4, 3) {};
		\node  (12) at (-3, 2) {};
		\node  (13) at (-2, 3) {};
		\node  (14) at (-2, 2) {};
		\node  (15) at (0, 2) {};
		\node  (16) at (0, 3) {};
		\node  (17) at (-1, 2) {};
		\node  (18) at (-6, 1) {};
		\node  (19) at (-6, 0) {};
		\node  (20) at (-5, 0) {};
		\node  (21) at (-4, 0) {};
		\node  (210) at (-3.6, 0.2) {$w_{15}$};
		\node  (22) at (-4, 1) {};
		\node  (23) at (-3, 0) {};
		\node  (24) at (-2, 0) {};
		\node  (25) at (-2, 1) {};
		\node  (26) at (0, 1) {};
		\node  (27) at (0, 0) {};
		\node  (28) at (-1, 0) {};
		\node  (29) at (1, 4) {};
		\node  (30) at (2, 4) {};
		\node  (31) at (2, 3) {};
		\node  (32) at (1, 2) {};
		\node  (33) at (2, 2) {};
		\node  (34) at (1, 0) {};
		\node  (35) at (2, 0) {};
		\node  (36) at (2, 1) {};
		\node  (37) at (-6, -1) {};
		\node  (38) at (-6, -2) {};
		\node  (39) at (-5, -2) {};
		\node  (40) at (-4, -2) {};
		\node  (41) at (-4, -1) {};
		\node  (42) at (-3, -2) {};
		\node  (43) at (-2, -2) {};
		\node  (44) at (-2, -1) {};
		\node  (45) at (0, -1) {};
		\node  (46) at (0, -2) {};
		\node  (47) at (1, -2) {};
		\node  (48) at (2, -2) {};
		\node  (49) at (2, -1) {};
		\node  (50) at (-1, -2) {};
		\node  (51) at (-6, -2.25) {$w_1$};
		\node  (52) at (-5, -2.25) {$w_2$};
		\node  (53) at (-4, -2.25) {$w_3$};
		\node  (54) at (-3, -2.25) {$w_4$};
		\node  (55) at (-2, -2.25) {$w_5$};
		\node  (56) at (-1, -2.25) {$w_6$};
		\node  (57) at (0, -2.25) {$w_7$};
		\node  (58) at (1, -2.25) {$w_8$};
		\node  (59) at (2, -2.25) {$w_1$};
		\node  (60) at (-6, 4.25) {$w_1$};
		\node  (61) at (-5, 4.25) {$w_2$};
		\node  (62) at (-4, 4.25) {$w_3$};
		\node  (63) at (-3, 4.25) {$w_4$};
		\node  (64) at (-2, 4.25) {$w_5$};
		\node  (65) at (-1, 4.25) {$w_6$};
		\node  (66) at (0, 4.25) {$w_7$};
		\node  (67) at (1, 4.25) {$w_8$};
		\node  (68) at (2, 4.25) {$w_1$};
		\node  (69) at (-6.25, 0) {$w_{13}$};
		\node  (70) at (-4.75, -0.25) {$w_{14}$};
		\node  (71) at (-3.75, -0.25) {};
		\node  (72) at (-2.75, -0.25) {$w_{16}$};
		\node  (73) at (-1.75, -0.25) {$w_{17}$};
		\node  (74) at (-0.75, -0.25) {$w_{18}$};
		\node  (75) at (0.25, -0.25) {$w_{19}$};
		\node  (76) at (1.25, -0.25) {$w_{20}$};
		\node  (77) at (2.25, -0.25) {$w_{13}$};
		\node  (78) at (-5.75, 1.75) {$w_{25}$};
		\node  (79) at (-4.75, 1.75) {$w_{26}$};
		\node  (80) at (-3.75, 1.75) {$w_{27}$};
		\node  (81) at (-2.75, 1.75) {$w_{28}$};
		\node  (82) at (-1.75, 1.75) {$w_{29}$};
		\node  (83) at (-0.75, 1.75) {$w_{30}$};
		\node  (84) at (0.25, 1.75) {$w_{31}$};
		\node  (85) at (1.25, 1.75) {$w_{32}$};
		\node  (86) at (2.25, 1.75) {$w_{25}$};
		\node  (87) at (-6.25, 3) {$w_{33}$};
		\node  (89) at (-4.5, 3) {$w_{34}$};
		\node  (91) at (-2.5, 3) {$w_{35}$};
		\node  (93) at (-0.5, 3) {$w_{36}$};
		\node  (95) at (1.5, 3) {$w_{33}$};
		\node  (96) at (-6.25, 1) {$w_{21}$};
		\node  (98) at (-4.5, 1) {$w_{22}$};
		\node  (100) at (-2.5, 1) {$w_{23}$};
		\node  (102) at (-0.5, 1) {$w_{24}$};
		\node  (104) at (1.5, 1) {$w_{21}$};
		\node  (105) at (-6.25, -1) {$w_9$};
		\node  (107) at (-4.25, -1) {$w_{10}$};
		\node  (109) at (-2.25, -1) {$w_{11}$};
		\node  (111) at (-0.25, -1) {$w_{12}$};
		\node  (113) at (1.75, -1) {$w_{9}$};
	
		\draw (0.center) to (6.center);
		\draw (6.center) to (27.center);
		\draw (27.center) to (19.center);
		\draw (19.center) to (0.center);
		\draw (8.center) to (15.center);
		\draw (1.center) to (21.center);
		\draw (2.center) to (24.center);
		\draw (7.center) to (3.center);
		\draw (18.center) to (4.center);
		\draw (20.center) to (5.center);
		\draw (23.center) to (16.center);
		\draw (28.center) to (26.center);
		\draw (19.center) to (38.center);
		\draw (38.center) to (48.center);
		\draw (48.center) to (30.center);
		\draw (30.center) to (6.center);
		\draw (16.center) to (29.center);
		\draw (15.center) to (33.center);
		\draw (27.center) to (35.center);
		\draw (39.center) to (23.center);
		\draw (20.center) to (37.center);
		\draw (21.center) to (40.center);
		\draw (24.center) to (43.center);
		\draw (27.center) to (46.center);
		\draw (42.center) to (28.center);
		\draw (26.center) to (31.center);
		\draw (36.center) to (50.center);
		\draw (47.center) to (49.center);
       }
\end{tikzpicture}
    }  
    \caption{Two isomorphic maps}
\end{figure}\label{eg}

From \cite{hubard:2012} we know that 
\begin{equation}\label{eqn-1}
   {\rm  Aut}\left(\frac{E}{K}
   \right)=\frac{\Nor_{\Aut(E)}(K)}{K}. 
\end{equation} This means, for any automorphism $\mathscr{h} \in \Aut(E)$, $\mathscr{h} \in \Aut(X)$ if and only if $\mathscr{h}$ normalizes $K$. In the following lemmas we find out the conditions on the matrix of $X$ for which an automorphism $\mathscr{h}\in \Nor(K)$. We also count the total number of such maps.

Now, to enumerate the number of 2-uniform maps we associate each symmetry in the tiling with a polynomial in $\ZZ_r[x]$ for some $r\in \mathbb{N}$ in the sense that if the entries of the HNF of a map satisfies the polynomial then that symmetry is present in its automorphism group.
Then to count the maps up to isomorphism we use the concept of isotropy group of a map as defined bellow.

\begin{definition}
For a toroidal map $X = E/K$, the {\em isotropy group of X} is defined as $\frac{\Aut(X)}{T_X}$, where $T_X =\Nor_H(K)/K = H/K$ where $H$ be the translation subgroup of $\Aut(E)$. We denote it by $\mathscr{I}(X)$. The isotropy group of the tiling $E$ is defined as $\Aut(E)/H$.
\end{definition}

Clearly, by (\ref{eqn-1}) $\mathscr{I}(X) \le \mathscr{I}(E)$. For example the isotropy group of the tiling of type $[3^1,6^1,3^1,6^1]$ is dihedral group of order 12, $D_6$. This is generated by $60^{\circ}$ rotation and a line reflection. The isotropy group of the map $M$ contains only reflection about a point. Thus $\mathscr{I}(M)$ is $\mathbb{Z}_2$. The number of distinct fundamental domains of this map keeping one vertex of the fundamental domain at $(0,0)$ is number of cosets of $\ZZ_2$ on $D_6$ which is 6.

\section{2-uniform toroidal maps}\label{proof}
In this section we proof our main theorems. First we will note the symmetries present in the 2-uniform tilings presented in Sec. \ref{sec:examples}.


Let $X_0$ be the minimal toroidal map, that is $E/H$ for some tiling $E$ of $\mathbb{R}^2$ and its translation group $H$.
Let $X$ be an $n$ sheeted cover of $X_0$, where $n\in\mathbb{N}$. 
Let $\tau$, $\sigma$, $\psi$ denote the function obtained by taking $180^{\circ}$, $60^{\circ}$ and $90^{\circ}$ rotation in $E$ with respect to $O$ respectively. For the  
tilings in Fig. \ref{Fig-1} and \ref{fig-2}
except $E_{16}$(\ref{E_16}) and $E_{23}$(\ref{E_23}), let $R_1$, $R_2$, $R_3$ denote the line passing through $O$ in the direction of the vectors $A$, $A-B$ and $B$. Let $R_4$, $R_5$ and $R_6$ denote the perpendicular bisectors of $R_1$, $R_2$, $R_3$ respectively. For Figures \ref{E_15}, \ref{E_16}, \ref{E_18}, \ref{E_19}, \ref{E_23} and \ref{E_27} let $R'_1$, $R'_2$ denote the line passing through $O$ in the direction of the vectors $A$ and perpendicular to it respectively. Let $R'_3$ denotes the line making an angle $45^{\circ}$ with $\overrightarrow{OA}$ for Fig. \ref{E_16} ,\ref{E_23} and \ref{E_27}. Let $R_4'$ be the line perpendicular to $R_3'$. Let $r_1$, $r_2$, $r_3$, $r_4$, $r_5$, $r_6$, $r'_1$, $r'_2$, $r'_3$ and $r'_4$ denote the functions obtained by taking reflections of $E$ about the lines $R_1$, $R_2$, $R_3$, $R_4$, $R_5$, $R_6$, $R'_1$, $R'_2$, $R'_3$ and $R'_4$ respectively. 

Before going into the details of the proofs let's observe the symmetries of the tilings and note which symmetries are responsible for determining level of symmetry of a map in other words presence of which symmetries in the automorphism group of the map changes the number of vertex orbits of the map.


Let us note the symmetries present in different tilings of the plane.
All the tilings have translations along the vectors $\overrightarrow{OA}$ and $\overrightarrow{OB}$ and point reflection with respect to $O$. Additionally,
For $i = 2,5,6,7,11,14,17,20,21,22,24~{\rm and}~26$ $E_i$ has the following bunch of symmetries and their compositions,
 $60^{\circ}$ rotation, $r_i$ for $i=1,2,\dots,6$. Here all of these functions are taking part in reducing orbits.
For $i=3,4,8,12,13,15,18,19$ $E_i$ has the following bunch of symmetries and their compositions,
reflection about horizontal and vertical lines through $O$. These functions are taking part in reducing orbits.
For $i=1$ and $25$ $E_i$ has only 
$60^{\circ}$ rotation and its compositions. This is also taking part in reducing orbits.
For $i=15,16,18,19,23$ and $27$ $E_i$ has the following bunch of symmetries and their compositions,
$90^{\circ}$ rotation, $r_i'$ for $i=1,2$. These are also taking part in reducing orbits. For $i=16,23,27$ together with these we also have $r_3'$ and $r_4'$. 
For $i=23$ $E_i$ has the following additional symmetries and their compositions,
Glide reflection whose component translation and reflection are not symmetries of the tiling. 


Observe that the isotropy group of the tilings $E_i$ is $D_6$ for $i \in \{2,5,6,7,11, 14,17,20,21,$ 
$22,24,26\}$, $D_4$ for $i \in \{16,27\}$, $\mathbb{Z}_4$ for $i=23$, $\mathbb{Z}_6$ for $i=1,25$ and $\mathbb{Z}_2 \times \mathbb{Z}_2$ for $i \in \{ 3,4,8,12,13,15,18,19 \}$. 
Now, for those tilings which have isotropy group $D_6$, the isotropy group of a corresponding map $X$ will be subgroup of $D_6$. So possible orders of $\mathscr{I}(X)$ will be $2,4,6 ~{\rm and}~ 12$. Order $1$ is not possible because point reflection is always present in $\Aut(X)$ and it has order $2$. 
There are $\sigma(n) = \sum_{d\mid n}d$ many distinct $n$-sheeted map exists. However, in this counting some cases counted more than once. If $\mathscr{I}(X)$ has order $2$ it is counted $6$ times (just by $60$ degree rotations). If $\mathscr{I}(X)$ has order $4$ then it is counted $3$ times. If $\mathscr{I}(X)$ has order $6$ it is counted $2$ times and if $\mathscr{I}(X)$ has order $12$ then it is counted only once. 
In general, if $\mathscr{I}(X)$ has order $d$ then it is counted $12/d$ times. This is caused because the choice of basis of the map is not unique, by applying different symmetries if the tiling on the basis we get other representation of the same map.
Similarly, for those tilings which have isotropy group $D_4$, the isotropy group of a corresponding map $X$ will be subgroup of $D_4$. So possible orders of $\mathscr{I}(X)$ will be $2,4 ~{\rm and}~ 8$. If $\mathscr{I}(X)$ has order $d$ then it is counted $8/d$ times.
And for those tilings which have isotropy group $\mathbb{Z}_4$ and $\mathbb{Z}_6$, the isotropy group of a corresponding map $X$ will be subgroup of $\mathbb{Z}_4$ and $\mathbb{Z}_6$ respectively. So possible orders of $\mathscr{I}(X)$ will be $2, 4$ and $2,6$. If $\mathscr{I}(X)$ has order $d$ then it is counted $4/d$ and $6/d$ times respectively. For the latter case order of $\mathscr{I}(X)$ cannot be $3$ since in that case the group will be generated by $120^{\circ}$ rotation and presence of $120^{\circ}$ rotation will ensure the presence of $60^{\circ}$ rotation.

Now, we are interested to know the symmetries present in Aut$(E)$ for which the vertex set of $E$ forms $2$-Nor$(K)$ orbits. 
Then we add all the counts corresponding to the possible isotropy groups of $X$. This gives the total number of $n$-sheeted $2$-orbital toroidal maps unique up to isomorphism. We denote this sum by denoted by $\Phi_i(v)$. 

First define,
$$\rho_{i,j}(n) := {\rm Number~of~solutions~of ~} x^2+ix+j\equiv 0 \pmod n.$$

Now, in the following subsections we discuss the classification of 2-uniform maps according to their types.

\subsection{Maps of type $[3^6;3^4,6^1]$($E_1$)}\label{first}

Let $X_0$ be the minimal toroidal map of type $[3^6;3^4,6^1]$. That is $X_0 = E_1/H_1$. Then any other map $X$ of this type will occur as a cover of $X_0$. 
Suppose $X$ be a $2$-orbital $n$-sheeted toroidal cover of $X_0$. Then, $X=\frac{E_1}{K_1}$ where $K_1\leq H_1$ with $H_1=\langle \gamma,\delta\rangle$ and $K_1=\langle w_1,\,w_2\rangle$. Since $X$ is a $2$-orbital map, $V(X)$ forms $2$-Aut$(X)$ orbits. Using  Result (\ref{eqn-1}),
we can say that $V(E_1)$ forms $2$- Nor$(K_1)$ orbits. So, we have to look for possible groups, $G = \Aut(E_1)$ such that $G\leq \Nor(K_1)$ and $V(E_1)$ forms $2$-$G$ orbits. 
Now, the translations $\gamma,\,\delta \in \Nor(K_1)$ since $\gamma^{-1}w_i\gamma=w_i$,  $\delta^{-1}w_i\delta=w_i$ $\in K_1$ for $i=1,2$.
We know that conjugation of a translation by rotation is again a translation by the rotated vector. So, $\tau^{-1}w_1\tau=-w_1 \in K_1$. Similarly, $\tau^{-1}w_2\tau\in K_1$. Hence, $\tau\in\Nor(K_1)$. Let $G_{1}=\langle \gamma,\,\delta,\,\tau\rangle$. Then, $V(E_1)$ forms $6$-$G_{1}$ orbits. Observe that $\rho\in $Aut$(E_1)$. 
Let $G_{2}=\langle \gamma,\,\delta,\,\tau,\,\rho \rangle$. 
Then, $V(E_1)$ forms $2$-$G_{2}$ orbits. Since $E_i$ has no reflection in its automorphism group, the only possible group giving 2 vertex orbits is $G_{2}$. So the isotropy group has order $6$. Hence, here we are counting every map exactly once.

\begin{lemma}\label{cm-1}
If $X$ be an $n$-sheeted toroidal cover of a given map $X_0$ of type $[3^6;3^4,6^1]$, then the total number of maps having $60^\circ$ rotation($\rho$) in the automorphism group of $X$ is equal to $$f_1(n)=\left\{
	\begin{array}{ll}
		0  & \mbox{, if } m_j\equiv1\pmod2 \mbox{ for some } j\in\{0,1,2,\dots,n_2\} \\
		\prod_{i=1}^{n_1}(k_i+1) & \mbox{, otherwise}
	\end{array}
\right.$$
where $n=2^{m_0} \cdot 3^{k_0} \cdot \prod_{i=1}^{n_1}p_i^{k_i} \cdot \prod_{j=1}^{n_2}q_j^{m_j}$ with $p_i$ and $q_j$ are primes such that $p_i\equiv1\pmod3$ for $i\in \{0,1,\dots,n_1\}$ and $q_j\equiv2\pmod3$ for $j\in \{0,1,\dots,n_2\}$.
\end{lemma}
\begin{proof}

Let $X$ be an $n$-sheeted toroidal cover of $X_0$ of type $[3^6;3^4,6^1]$ with associated matrix $\begin{bmatrix} a & 0\\ b & d\end{bmatrix}$. We have taken this matrix representation with respect to basis vectors of $X_0$. So, $ad=n$. Then $X=\frac{E}{\mathcal{K}}$ with $\mathcal{K}=\langle w_1,\,w_2\rangle$ where $w_1:z\mapsto z+(aA+bB)$ and $w_2: z\mapsto z+dB$ (See Fig. \ref{Fig-1}). 
Let $\rho$ denotes the transformation obtained by taking $60^{\circ}$ rotation on $E_1$.
$\rho$ sends $A$ to $B$ and $B$ to $B-A$.

We look for the conditions on $a,\,b$ and $d$ so that $\rho\in\Nor(\mathcal{K})$. For that, it is enough to check if $\rho^{-1}w_1\rho$, $\rho^{-1}w_2\rho$ belong $\mathcal{K}$. 
We know that conjugation of a translation by a rotation is again a translation by the rotated vector. So, $\rho^{-1}w_i\rho$  is translation by the vector $\rho \circ w_i(0).$  Now, $\rho \circ w_1(0)=\rho(aA+bB)=a\rho(A)+b\rho(B)=aB+b(B-A)=-bA+(a+b)B$ and  $\rho \circ w_2(0)=\rho(cA+dB)=c\rho(A)+d\rho(B)=cB+d(B-A)=-dA+(c+d)B$. Now, $\rho(w_1)$ and $\rho(w_2)$ belong to the lattice of $\mathcal{K}$ provided that there exists integers $m_1,\,m_2,\,m_3,\, m_4$ such that $-bA+(a+b)B=m_1(aA+bB)+m_2dB$ and $-dA+dB=m_3(aA+bB)+m_4dB.$ Since, $A$ and $B$ are linearly independent, we have a system of linear equations, $m_1a=-b$, $m_1b+m_2d=a+b$, $m_3a=-d$ and $m_3b+m_4d=d$. Solving these equations, we get, $m_1=-\frac{b}{a}$, $m_2=\frac{a^2+ab+b^2}{ad}$, $m_3=-\frac{d}{a}$ and $m_4=\frac{a+b}{a}$. Since $m_1,\,m_2,\,m_3$ and $m_4$ are integers, we must have $a\mid b$, $a\mid d$ and $ad\mid(a^2+ab+b^2)$. Therefore, $\rho \in \Aut(X)$ if and only if {\rm(i)} $a\mid b$, {\rm(ii)} $a\mid d$ and {\rm(iii)} $ad\mid(a^2+ab+b^2)$.

Now, $a\mid b$ and $a \mid d$ implies $b=ax$ and $d=ay$ where $0\leq x<y$. Using the last condition, we get $a.ay\mid (a^2+a.ax+ax.ax)$, or, $y\mid (1+x+x^2)$ where $y=\frac{d^2}{n}$. This implies $ 1+x+x^2$ has a solution in $\mathbb{Z}_{\frac{d^2}{n}}$. So, the total number of maps is given by the number of solutions of the polynomial $1+x+x^2$ in $\mathbb{Z}_{\frac{d^2}{n}}$ satisfying $n\mid d^2$ for every divisor $d$ of $n$. Hence, we have 
$$f_1(n):=\sum_{d \mid n,\,n \mid d^2}\rho_{1,1}\left(\frac{d^2}{n}\right).$$
From \cite[Chapter 8, Theorem-122]{HW1979} we have the following,
\begin{result}\label{result--1}
The number of roots of $f(x) \equiv 0 \pmod n$ is the product of the number of roots of separate congruences $f(x) \equiv 0 \pmod{p_i^{t_i}}$ for $i=1,2,3,\dots,k$ where $n=p_1^{t_1}p_2^{t_2}\cdots p_k^{t_k}$.
\end{result}
Applying the above result we get $\rho_{1,1}$ is a multiplicative function. As a consequence of this we have $f_1$ is also multiplicative. Thus it is enough to find values of $f_1$ for prime powers. In the following $\lceil x \rceil$ denotes the smallest integer grater than or equals to $x$.

$f_1(p^k)=\sum_{\substack{d \mid p^k \\p^k \mid d^2}}\rho_{1,1}\left(\frac{d^2}{p^k}\right) = 
\sum_{\substack{p^l \mid p^k \\p^k \mid p^{2l}}}\rho_{1,1}\left(\frac{p^{2l}}{p^k}\right)=
\sum_{\frac{k}{2}\le l \le k}\rho_{1,1}\left(p^{2l-k}\right)=
\sum_{l=\lceil \frac{k}{2} \rceil}^k\rho_{1,1}\left(p^{2l-k}\right).$
From theory of quadratic residues \cite{HW1979} we have $1+x+x^2 \equiv 0 \pmod {p^i}$ has 2 distinct solutions when $p \equiv 1 \pmod  3$ provided $i>0$. If $p \equiv 2 \pmod 3$ and odd then $1+x+x^2 \equiv 0 \pmod {p^i}$ has no solution provided $i>0$. If $i=0$ then it has unique solution in both the cases.
Putting these values in the above formula we have,
$$f_1(p^k) =\left\{
	\begin{array}{ll}
		\sum_{l= \frac{k}{2} }^k\rho_{1,1}\left(p^{2l-k}\right)  & \mbox{, if  2 $\mid$ k} \\
		\sum_{l= \frac{k+1}{2} }^k\rho_{1,1}\left(p^{2l-k}\right) & \mbox{, if  2 $\nmid$ k}
	\end{array}
          \right. =\left\{
	\begin{array}{lll}
		k+1  &, p \equiv 1 \pmod  3 \\
		0 &, p \equiv 2 \pmod  3 ~{\rm and}~ 2 \nmid k \\
		1 &, p \equiv 2 \pmod  3 ~{\rm and}~ 2 \mid k.
	\end{array}
          \right.$$
          
Now,
observe that $\rho_{1,1}(2^k)=0$ for all $k\in \mathbb{N}$ since $1+x+x^2$ is odd for every $x\in \mathbb{N}\cup\{0\}$ and $2^k\nmid (1+x+x^2)$.
$f_1(3^k)=1$ for all $k\in\mathbb{N}$.
$\rho_{1,1}(3)=1$ since $1$ is the only solution of the equation $1+x+x^2$ in $\mathbb{Z}_3$.
We want to show that $\rho_{1,1}(3^k)=0$ for all $k\geq 2$. When $x\equiv 0 \pmod 3$ then $1+x+x^2\equiv 1 \pmod 3$. When $x=3l+1$ where $l\in\mathbb{N}$, $1+x+x^2=3(1+3l+3l^2)$. When $x\equiv 2 \pmod 3$ then $1+x+x^2\equiv 1 \pmod 3$. So, $3^k\nmid(1+x+x^2)$ for all $x\in \mathbb{N}$. Hence, $\rho_{1,1}(3^k)=0$. Therefore we have,
$$f_1(2^k) =\left\{
	\begin{array}{ll}
		\sum_{l= \frac{k}{2} }^k\rho_{1,1}\left(2^{2l-k}\right) = 1 & \mbox{, if  2 $\mid$ k} \\
		\sum_{l= \frac{k+1}{2} }^k\rho_{1,1}\left(2^{2l-k}\right) = 0  & \mbox{, if  2 $\nmid$ k,}
	\end{array}
          \right. 
f_1(3^k) =\left\{
	\begin{array}{ll}
		\sum_{l= \frac{k}{2} }^k\rho_{1,1}\left(3^{2l-k}\right) = 1 & \mbox{, if  2 $\mid$ k} \\
		\sum_{l= \frac{k+1}{2} }^k\rho_{1,1}\left(3^{2l-k}\right) = 1  & \mbox{, if  2 $\nmid$ k.}
	\end{array}
          \right. $$

Now, by fundamental theorem of arithmetic any integer $n$ can be written as  $n=2^{m_0} \cdot 3^{k_0} \cdot \prod_{i=1}^{n_1}p_i^{k_i} \cdot \prod_{j=1}^{n_2}q_j^{m_j}$ with $p_i$ and $q_j$ are primes such that $p_i\equiv1\pmod3$ for $i\in \{0,1,\dots,n_1\}$ and $q_j\equiv2\pmod3$ for $j\in \{0,1,\dots,n_2\}$, using multiplicative property of $f_1$ the assertion of Lemma \ref{cm-1} follows.
\end{proof}

\begin{remark}
Lemma \ref{cm-1} is valid for other type of maps also. That is for a type $Z$ the number of distinct maps having $\rho$ in its automorphism group is $f_1(n)$ where $n$ is the number of sheets.
\end{remark}

\begin{remark}
$f_1$ is a special kind of divisor function. It counts number of divisors of $n>1$ which are congruent to 1 modulo 3.
\end{remark}

Therefore, total number of $2$-orbital $n$ sheeted covers up to isomorphism is $f_1(n)$. Note that there are 12 vertices in $X_0$. Thus the number of 2-uniform maps up to isomorphism with $v$ vertices is,
\begin{center}
    $\Phi_{1}(v) = \left\{
	\begin{array}{ll}
		0  & \mbox{, if } 12 \nmid v \\
		f_1(\frac{v}{12}) & \mbox{, otherwise.}
	\end{array}
      \right.$ 
\end{center}

\subsection{Maps of type $[3^6;3^4,6^1](E_2),[3^6;3^2,4^1,3^1,4^1],[3^6;3^2,6^2],[3^{3},4^2;3^1,4^1,6^1,4^1],$\\$[3^1, 4^1, 6^1, 4^1;3^2,4^1,3^1,4^1]$}


Let $X_0$ be the minimal toroidal map of type any one of the above. That is $X_0 = E_i/H_i$ for $i=2,5,7,11$ or14(see Fig. \ref{Fig-1}).
Let $X$ be an $n$ sheeted $2$-orbital cover of $X_0$. Then, $X=\frac{E_i}{K_i}$ where $K_i\leq H_i$ with $H_i=\langle \gamma,\delta\rangle$ and
$K_i=\langle w_1,\,w_2\rangle$. Since $X$ is a $2$-orbital map, $V(X)$ forms $2$-Aut$(X)$ orbits. Using Result (\ref{eqn-1}),
we can say that $V(E_i)$ forms $2$-Nor$(K_i)$ orbits. 
Now $V(E_i)$ forms $4$ and $6$ orbits for $i=2,5,7$ and $i=11,14$ respectively when only $\tau$ and translations present in  Nor$(K_i)$. If $\mathscr{I}(X)$ contains $\tau$ and $r_j$ for some $j\in \{1,2,\dots,6\}$ then $V(E_i)$ has $3$ Nor($K_i$) orbits. If $\mathscr{I}(X)$ contains $\rho$ then $V(E_i)$ has $2$ Nor($K_i$) orbits. And if $\mathscr{I}(X)$ is equals to isotropy group of the tiling then also it has $2$ orbits. 

By Lemma \ref{cm-1} we can say that the number of distinct maps having $\rho$ in its automorphism group is $f_1(n)$. Now we have to find the number of distinct maps having its isotropy group equals to that of the tiling. If the isotropy group of the map contain $r_1$ and $r_2$ then it will be equal to the isotropy group of the tiling. The following lemma gives number of distinct maps with this property.

\begin{lemma}\label{lem-f5}
If $X$ be a $n$ sheeted toroidal cover of $X_0$, then the total number of maps having $r_1$ and $r_2$ in the automorphism group of X is equal to 
$$f_5(n)=\left\{
	\begin{array}{ll}
		0  & \mbox{, if } k_i\equiv1\pmod2 \mbox{ for some } i\in\{0,1,2,\dots,n_1\}\\
		1 & \mbox{, otherwise}
	\end{array}
\right.$$
where $n=2^{k_0}3^{m_0}p_1^{k_1}p_2^{k_2}\dots p_{n_1}^{k_{n_1}}$ where $p_i$ is any prime other than $3$ for $i\in \{0,1,\dots,n_1\}$.
\end{lemma}

\begin{proof}
Let $X$ be a $n$ sheeted toroidal cover of $X_0$ with HNF of the associated matrix is $\begin{bmatrix}a&0\\b&d\end{bmatrix}$. Let $X$ has $r_1$ and $r_2$ in its automorphism group.
The matrices of $r_1$ and $r_2$ are $\begin{bmatrix} 1&1\\0&-1\end{bmatrix}$ and $\begin{bmatrix} -1&1\\1&0\end{bmatrix}$ respectively. Proceeding similar way as in the proof of Lemma \ref{cm-1} we get 
$r_1$ and $r_2$ belong to Aut$(X)$ if and only if $a\mid b$, $a\mid d$, $ad\mid (b^2+2ab)$ and $ad\mid(b^2-a^2)$. Let $b=ax$ and $d=ay$.

Using these, we get $y\mid (x^2+2x)$ and $y\mid (1-x^2)$ where $y=\frac{d^2}{n}$ and $0\leq x<y$. This implies that the polynomials $x^2+x$ and $1-x^2$ have a common solution in $\mathbb{Z}_{\frac{d^2}{n}}$. So, the total number of maps is given by the number of solutions of gcd$(x^2+2x,\,1-x^2)$ in $\mathbb{Z}_{\frac{d^2}{n}}$ satisfying $n\mid d^2$ for every divisor $d$ of $n$. Hence, we have,
$$f_5(n):=\sum_{\substack{d \mid n,~n \mid d^2}}\rho_5\left(\frac{d^2}{n}\right)$$ 
where $\rho_5(n):=\#\{x\in \mathbb{Z}_n:\, x^2+2x=0\,{\rm and}\,1-x^2=0\}$.

By Result \ref{result--1} we have $\rho_5$ is multiplicative. Thus $f_5$ is so. Similarly as in Lemma \ref{cm-1} $f_5(p^k) = \sum_{l=\lceil \frac{k}{2} \rceil}^k\rho_{5}\left(p^{2l-k}\right).$

Now, observe that $x=0$ is the solution for the congruences $x^2+2x\equiv0\pmod1$ and $1-x^2\equiv0\pmod1$. Thus $\rho_{5}(1)=1$. since $gcd\{(x^2+2x),\,(x^2-1)\}=x-1$ in $\mathbb{Z}_3[x]$, $x=1$ is the only solution in $\mathbb{Z}_3$. Thus, $\rho_{5}(3)=1$.
Since, gcd$\{(x^2+2x),\,(x^2-1)\}=1$ in $\mathbb{Z}_{3^k}[x]$ where $k\in \mathbb{N}\setminus\{1\}$, we have $\rho_5(3^k)=0$ for all $k \neq 1$.
Therefore, $f_5(3^k)=1$ for all $k\in\mathbb{N}$.
If $p$ in any prime other than $3$, then gcd$\{(x^2+2x),\,(x^2-1)\}=1$ in $\mathbb{Z}_{p^k}[x]$. Thus we have $\rho_5(p^k)=0$ for all $k\in \mathbb{N}$. Hence $f_5(p^k) = 1$ for all $k\in \mathbb{N}$.

Now using multiplicative property of $f_5$ the assertion follows.
\end{proof}

Now, depending on the order of the isotropy group those maps having only $\rho$ in its automorphism group are counted twice in the list of distinct maps and those maps having full isotropy group are counted once.
Hence number of $n$-sheeted $2$-uniform maps up to isomorphism there are $\frac{f_1(n)-f_5(n)}{2}+ f_5(n) = \frac{f_1(n)+ f_5(n)}{2} $.
Let there are $v_0$ vertices in $X_0$. Thus the number of 2-uniform maps up to isomorphism with $v$ vertices is $\Phi_{\ell}$  given by,
$\Phi_{\ell}(v) = \left\{
	\begin{array}{ll}
		0  & \mbox{, if } v_0 \nmid v \\
		\frac{1}{2}[f_1(\frac{v}{v_0})+f_5(\frac{v}{v_0})] & \mbox{, otherwise}
	\end{array}
      \right.$ for  $(\ell,\,v_0) = (2,8),(5,7),(7,7),(11,12),(14,12)$.

\subsection{Maps of type $[3^6;3^3,4^{2}] (E_3 ~\&~ E_4), [3^2,6^2;3^4,6^1], [4^4;3^3,4^2], (E_{12} \& E_{13}),$\\ $[3^2,6^2;3^1, 6^1, 3^1, 6^1]$}


Let $X$ be an $n$ sheeted $2$-orbital cover of $X_0$ obtained as an orbit space of $E_i$ where $i=3,\,4,\,8,\,12,\,13,\,15$ (see Fig. \ref{Fig-1}). 
In these cases every map is $2$-orbital. So, it is enough to find total number of $n$-sheeted maps up to isomorphism.
Isotropy group of these tilings is $\ZZ_2 \times \ZZ_2$. So, possible orders $\mathscr{I}$(X) will be $2$ and $4$ according to the presence of the line reflection along the vector $\overrightarrow{OA}$. If it is $4$ then we count it once in the collection of distinct maps and if it is $2$ then we count it twice. Let $C(n)$ and $D(n)$ denoted number of maps up to isomorphism having order of isotropy group $2$ and $4$ respectively. Let $h$ denotes the reflection.
Then $2C(n) + D(n) = \sigma(n)$. Now, $D(n)$ is also number of maps having $h$ in its automorphism group.
Now, depending on the different tilings the matrix of $h$ is different. We discuss them case by case in the following.

\subsubsection{Maps of type $[3^2,6^2;3^4,6^1], [4^4;3^3,4^2](E_{12} \& E_{13})$}

Here the matrix of $h$ will be $\begin{bmatrix} 1&1\\0&-1\end{bmatrix}$. The following lemma gives the number of maps having $h$ in its isotropy group.

\begin{lemma}\label{lem--3}
If $X$ be a $n$ sheeted toroidal cover of $X_0$, then the total number of maps having $h$ in the automorphism group of $X$ is equal to,
\[ f_3(n)=\left\{
	\begin{array}{ll}
		\prod_{i=1}^{n_1}(k_i+1)  & \mbox{, if } k_0=0 \\
		(2k_0-1)\prod_{i=1}^{n_1}(k_i+1) & \mbox{, otherwise}
	\end{array}
\right. \]
where $n=2^{k_0}p_1^{k_1}p_2^{k_2}\dots p_{n_1}^{k_{n_1}}$ where $p_i$ is any odd prime for $i\in \{0,1,\dots,n_1\}$.

\end{lemma}
\begin{proof}
Let the associated matrix of $X$ be $\begin{bmatrix} a & 0\\ b & d\end{bmatrix}$. 
Now proceeding in similar way as in Lemma \ref{cm-1} one can see that 
$h\in$ $\Nor(K)$ if and only if {\rm (i)} $a\mid b$, {\rm (ii)} $a\mid d$ and {\rm (iii)} $ad\mid (b^2+2ab)$, 
Now, let number of distinct maps having $h$ in its isotropy group is $f_3(n)$. Then, 
$f_3(n)=\sum_{d \mid n,~n \mid d^2}\rho_{2,0}\left(\frac{d^2}{n}\right).$
By similar reasons $f_3$ is also multiplicative.
Observe that,
$\rho_{2,0}(1)=1$ since $x=0$ is the solution for $x^2+2x\equiv0\pmod1$.
$\rho_{2,0}(2)=1$ since $0$ is the only solution of $0$ in $\mathbb{Z}_2$.
$\rho_{2,0}(4)=2$ since $0$ and $2$ are the solutions of $x^2+2x=0$ in $\mathbb{Z}_4$.
Now, $x^2+2x\equiv 0\pmod{2^k}\implies (x+1)^2\equiv1\pmod{2^k}$. Putting $y=x+1$, we get the congruence $y^2\equiv1\pmod{2^k}$. We know that the congruence $y^2\equiv1\pmod{2^k}$ has exactly four incongruent solutions. Hence, $\rho_{2,0}(2^k)=4$.
The congruence $x^2+2x\equiv0\pmod{p}$ has solution if and only if $y^2\equiv4\pmod{p}$ has a solution. Since, $p$ is odd we have $\left(\frac{4}{p}\right)=1$ that is 4 is a quadratic residue $\mod p$. Hence, $x^2+2x\equiv0\pmod{p}$ has exactly two solutions for all odd prime $p$. This implies $x^2+2x\equiv0\pmod{p^k}$ has $2$ solutions for all $k\geq1$. So, $\rho(p^k)=2$ for all $k\geq1$.
Therefore,
$$f_3(2^k) = 2k-1 ~ {\rm and} ~ 
f_3(p^k) = k+1
.$$
Therefore by multiplicative property of $f_3$ the assertion follows.
\end{proof}
Hence, $D(n) = f_3(n)$. Thus, number of 2-uniform maps up to isomorphism if given by $C(n) + D(n) = \frac{\sigma(n) + f_3(n)}{2}$.
Let there are $v_0$ vertices in $X_0$. Thus the number of 2-uniform maps up to isomorphism with $v$ vertices is $\Phi_{\ell}$  given by,
$\Phi_{\ell}(v) =\left\{
	\begin{array}{ll}
		0  & \mbox{, if } v_0 \nmid v \\
		\frac{1}{2}[\sigma (\frac{v}{v_0})+f_3(\frac{v}{v_0})] & \mbox{, otherwise}  
	\end{array}
      \right.$  for $(\ell, v_0)=(8,4),(12,3),(13,4)$.

\subsubsection{Maps of type $[3^6;3^3,4^{2}](E_3)$}       
Here the matrix of $h$ will be $\begin{bmatrix} 1&3\\0&-1\end{bmatrix}$. The following lemma gives the number of maps having $h$ in its isotropy group.
Let $X$ be a map of this type with associated matrix $\begin{bmatrix} a&0\\b&d\end{bmatrix}$. Then by similar type of calculation as in Lemma \ref{cm-1} we get $h$ will present in Aut($X$) if and only if $a \mid 3b$, $a\mid 3d$ and $ad\mid 3b^2+2ab$.
Let $g_2(n):=\#\left\{\begin{bmatrix} a&0\\b&d\end{bmatrix} : a\mid 3b, a\mid 3d ~{\rm and}~ ad \mid 3b^2+2bd\right\}$. Then $g_2$ is a multiplicative function. Therefore it is enough to find the values $g_2(p^k)$ for primes $p$.
For odd primes except 3, $g_2(p^k) = \sum_{i=\lceil \frac{k}{2}\rceil}^k \#\{x:x(3x+2)\equiv0 \pmod{p^{2i-k}}\} = \left\{
	\begin{array}{ll}
		k+2  & \mbox{, if } 2 \mid k \\
		k+1 & \mbox{, otherwise.}  
	\end{array}
      \right. $
For even prime, $g_2(2^k) =  \left\{
	\begin{array}{ll}
		k+2  & \mbox{, if } 2 \mid k \\
		k & \mbox{, otherwise.}  
	\end{array}
      \right. $
Finally, 
$g_2(3^k) =  k+2 $ for all $k\in \mathbb{N}$.
Hence, $D(n) = g_2(n)$. Thus, number of 2-uniform maps up to isomorphism if given by $C(n) + D(n) = \frac{\sigma(n) + g_2(n)}{2}$. Therefore,
$\Phi_{3}(v) =\left\{
	\begin{array}{ll}
		0  & \mbox{, if } 4 \nmid v \\
		\frac{1}{2}[\sigma (\frac{v}{4})+g_2(\frac{v}{4})] & \mbox{, otherwise.}  
	\end{array}
      \right.$     

\subsubsection{Maps of type $[3^6;3^3,4^{2}](E_4)$}       
Here the matrix of $h$ will be $\begin{bmatrix} 1&2\\0&-1\end{bmatrix}$. The following lemma gives the number of maps having $h$ in its isotropy group.
Let $X$ be a map of this type with associated matrix $\begin{bmatrix} a&0\\b&d\end{bmatrix}$. Then by similar type of calculation as in Lemma \ref{cm-1} we get $h$ will present in Aut($X$) if and only if $a \mid 2b$ and $a\mid 2d$.
Let $g_1(n):=\#\left\{\begin{bmatrix} a&0\\b&d\end{bmatrix} : a\mid 2b, a\mid 2d\right\}$. Then $g_1$ is a multiplicative function. Therefore it is enough to find the values $g_1(p^k)$ for primes $p$. 
For odd primes $g_1(p^k) = \sum_{i=\lceil \frac{k}{2}\rceil}^k p^{2i-k} = \left\{
	\begin{array}{ll}
		\frac{p^{k+2}-1}{p^2-1}  & \mbox{, if } 2 \mid k \\
		p\cdot\frac{p^{k+1}-1}{p^2-1} & \mbox{, otherwise.}  
	\end{array}
      \right. $
For even prime we have $g_1(2^k) = \sum_{i=\lceil \frac{k-1}{2}\rceil}^k 2^{2i+1-k} =\left\{
	\begin{array}{ll}
		\frac{2^{k+3}-1}{3}  & \mbox{, if } 2 \nmid k \\
		2\cdot\frac{2^{k+2}-1}{3} & \mbox{, otherwise.}  
	\end{array}
      \right. $
Hence, $D(n) = g_1(n)$. Thus, number of 2-uniform maps up to isomorphism if given by $C(n) + D(n) = \frac{\sigma(n) + g_1(n)}{2}$. Therefore,
$\Phi_{4}(v) =\left\{
	\begin{array}{ll}
		0  & \mbox{, if } 3 \nmid v \\
		\frac{1}{2}[\sigma (\frac{v}{3})+g_1(\frac{v}{3})] & \mbox{, otherwise.}  
	\end{array}
      \right.$

\subsubsection{Maps of type $[3^2,6^2;3^1, 6^1, 3^1, 6^1]$}\label{r_1}   
Here the matrix of $h$ will be $\begin{bmatrix} 1&0\\0&-1\end{bmatrix}$. 
Let $X$ be a map of this type with associated matrix $\begin{bmatrix} a&0\\b&d\end{bmatrix}$. Then by similar type of calculation as in Lemma \ref{cm-1} we get $h$ will present in Aut($X$) if and only if $d \mid 2b$. Therefore the number of such maps is $g(n) := \sum_{\substack{d \mid n,~2 \mid d}}2 + \sum_{\substack{d \mid n,~2 \nmid d}}1.$
Hence, $D(n) = g(n)$. Thus, number of 2-uniform maps up to isomorphism if given by $C(n) + D(n) = \frac{\sigma(n) + g(n)}{2}$. Therefore,
$\Phi_{15}(v) =\left\{
	\begin{array}{ll}
		0  & \mbox{, if } 3 \nmid v \\
		\frac{1}{2}[\sigma (\frac{v}{3})+g(\frac{v}{3})] & \mbox{, otherwise.}  
	\end{array}
      \right.$

\subsection{Maps of type $ [3^3, 4^2;3^2,4^1,3^1,4^1] (E_9 ~\& ~E_{10})$}


Let $X$ be an $n$ sheeted $2$-orbital cover of $X_0$ obtained as an orbit space of $E_i$ where $i=9,\,10$ (see Fig. \ref{Fig-1}).
Then, $X=\frac{E_i}{K_i}$ where $K_i\leq H_i$ with $H_i=\langle \gamma,\delta\rangle$ and
$K_i=\langle w_1,\,w_2\rangle$. Since $X$ is a $2$-orbital map, $V(X)$ forms $2$-Aut$(X)$ orbits. Using Result (\ref{eqn-1}),
we can say that $V(E_i)$ forms $2$- Nor$(K_i)$ orbits. 
In these two tilings we have two additional symmetries. We denote them by $\mathscr{G}_1$ and $\mathscr{G}_2$ and define them as follows.
$\mathscr{G}_1$ is the function obtained by rotating  the vertical stack of two squares by $90^{\circ}$ followed by a translation of that stack to its adjacent horizontal array of two squares in $E_{9}$.
$\mathscr{G}_2$ is the function obtained by reflecting $E_{10}$ about some line and then translating it by distance $|A|/2$ or $|B|/2$ along $\overrightarrow{OA}$ or $\overrightarrow{OB}$ respectively.
The maps will be $2$-uniform if and only if $\mathscr{G}_1$ and $\mathscr{G}_2$ present in the corresponding automorphism group. The following lemma gives the number of distinct maps having $\mathscr{G}_i$ in its automorphism group.

\begin{lemma}\label{glide}
Let $X$ be a $n$ sheeted toroidal map represented by the matrix $\begin{bmatrix}
a & 0 \\ b & d \end{bmatrix}$. Then $\mathscr{G}_i$ will be present in $\Aut(X)$ if and only if $d \mid 2b$ and number of distinct $n$-sheeted maps having glide in its automorphism group is 
$$g(n) := \sum_{\substack{d \mid n,~2 \mid d}}2 + \sum_{\substack{d \mid n,~2 \nmid d}}1.$$
\end{lemma}

\begin{proof}
Corresponding matrix for $\mathscr{G}_i$ is $\begin{bmatrix}
-1 & 0 \\ 0 & 1 \end{bmatrix}$. Now proceed similarly as Lemma \ref{cm-1}.
\end{proof}

\begin{remark}\label{remark1}
For $n \in \mathbb{N}$ let $k$ be the largest integer such that $2^k \mid n$. Then $g(n) = 2\tau(n) - \tau(n/2^k)$.
\end{remark}

Since in these tilings there does not exists any other symmetries like rotation or reflection about lines so this $g(n)$ is giving number of maps up to isomorphism. Thus number of $n$-sheeted $2$-uniform maps up to isomorphism is given by $g(n)$.

Let there are $v_0$ vertices in $X_0$. Thus the number of 2-uniform maps up to isomorphism with $v$ vertices is $\Phi_{\ell}$  given by,
$\Phi_{\ell}(v) =\left\{
	\begin{array}{ll}
		0  & \mbox{, if } v_0 \nmid v \\
		g (\frac{v}{v_0}) & \mbox{, otherwise}  
	\end{array}
      \right.$        for $(\ell,v_0) = (9,12),(10,8)$.

\subsection{Maps of type $ [3^1,4^2, 6^1; 3^1, 6^1, 3^1, 6^1](E_{18} ~\& ~E_{19})$}


Let $X$ be an $n$ sheeted $2$-orbital cover of $X_0$ obtained as an orbit space of $E_i$. 
$i=18,\,19$ (see Fig. \ref{Fig-1}). 
Then, $X=\frac{E_i}{K_i}$ where $K_i\leq H_i$ with $H_i=\langle \gamma,\delta\rangle$ and
$K_i=\langle w_1,\,w_2\rangle$. Since $X$ is a $2$-orbital map, $V(X)$ forms $2$-Aut$(X)$ orbits. Using Result (\ref{eqn-1}),
we can say that $V(E_i)$ forms $2$- Nor$(K_i)$ orbits.
It can be observed that these type of maps will be $2$-uniform if and only if $r_1'$ or $r_2'$ present in the automorphism group. We are counting them exactly once since their isotropy group is same as that of the tiling. 
Now, depending on the different tilings the matrix of $r_1'$ is different. We discuss them case by case in the following.

\subsubsection{Maps coming from  $E_{18}$}
Here the matrix of $r_1'$ will be $\begin{bmatrix} 1&0\\0&-1\end{bmatrix}$. Similarly as for maps of type $[3^2,6^2;3^1,6^1,3^1,6^1]$ we have number of maps having $r_1'$ in its automorphism group is $g(n) := \sum_{\substack{d \mid n,~2 \mid d}}2 + \sum_{\substack{d \mid n,~2 \nmid d}}1.$
Note that there are $5$ vertices in $X_0$. Thus the number of 2-uniform maps up to isomorphism with $v$ vertices is $\Phi_{18}$  given by,
$\Phi_{18}(v) =\left\{
	\begin{array}{ll}
		0  & \mbox{, if } 5 \nmid v \\
		g(\frac{v}{v_0}) & \mbox{, otherwise}  
	\end{array}
      \right.$.  

\subsubsection{Maps coming from  $E_{19}$}
Here the matrix of $r_1'$ will be $\begin{bmatrix} 1&1\\0&-1\end{bmatrix}$.
Hence by Lemma \ref{lem--3}  number of $n$-sheeted $2$-uniform maps up to isomorphism is given by $f_3(n)$. 

Note that there are $5$ vertices in $X_0$. Thus the number of 2-uniform maps up to isomorphism with $v$ vertices is $\Phi_{19}$  given by,
$\Phi_{19}(v) =\left\{
	\begin{array}{ll}
		0  & \mbox{, if } 5 \nmid v \\
		f_3(\frac{v}{v_0}) & \mbox{, otherwise}  
	\end{array}
      \right.$.

\subsection{Maps of type $ [3^{6}; 3^2, 4^1,12^1], [3^1,4^2,6^1; 3^1, 4^1, 6^1, 4^1], [4^1,6^1,12^1;3^1, 4^1, 6^1, 4^1]$}


Let $X$ be an $n$ sheeted $2$-orbital cover of $X_0$ obtained as an orbit space of $E_i$. 
$i=6,\,17,\,20$ (see Fig. \ref{Fig-1}). 
Then, $X=\frac{E_i}{K_i}$ where $K_i$ is as above. Similarly,
we can say that $V(E_i)$ forms $2$- Nor$(K_i)$ orbits.
For these types of maps one can observe that $V(E_i)$ has 2-Nor($K_i$) orbits if $\mathscr{I}(X) = \mathscr{I}(E_i)$. Hence by Lemma \ref{lem-f5} number of $2$-uniform maps up to isomorphism is $f_5(n)$.

Let there are $v_0$ vertices in $X_0$. Thus the number of 2-uniform maps up to isomorphism with $v$ vertices is $\Phi_{\ell}$  given by,
$\Phi_{\ell}(v) =\left\{
	\begin{array}{ll}
		0  & \mbox{, if } v_0 \nmid v \\
		f_5 (\frac{v}{v_0}) & \mbox{, otherwise}
	\end{array}
      \right.$         for $(\ell,v_0)=(6,14),(17,18),(20,18)$.

\subsection{Maps of type $[3^1,12^2;3^1,4^1,3^1,12^1]$}


Let $X$ be an $n$ sheeted $2$-orbital cover of $X_0$ obtained as an orbit space of $E_{16}$ (see Fig. \ref{Fig-1}).
Then, $X=\frac{E_{16}}{K_{16}}$ where $K_{16}\leq H_{16}$ with $H_{16}=\langle \gamma,\delta\rangle$ and
$K_{16}=\langle w_1,\,w_2\rangle$. Since $X$ is a $2$-orbital map, $V(X)$ forms $2$-Aut$(X)$ orbits. Using Result (\ref{eqn-1}),
we can say that $V(E_{16})$ forms $2$- Nor$(K_{16})$ orbits.
Now $V(E_i)$ forms $4$-Nor$(K_i)$  orbits when only $\tau$ and translations present in  Nor$(K_i)$. If $\mathscr{I}(X)$ contains $\tau$ and $r'_j$ for some $j\in \{1,2,\dots,4\}$ then $V(E_i)$ has $3$ Nor($K_i$) orbits. If $\mathscr{I}(X)$ contains $\psi$ then $V(E_i)$ has $2$ Nor($K_i$) orbits. And if $\mathscr{I}(X)$ is equals to isotropy group of the tiling then also it has $2$ orbits. 

The following lemmas gives the number of distinct maps having only $\psi$ in its isotropy groups and number of distinct maps having $\mathscr{I}(X) = \mathscr{I}(E_i)$.

\begin{lemma}\label{lem-f2}
If $X$ be an $n$-sheeted toroidal cover of a given map $X_0$, then the total number of maps having $\psi$ in the automorphism group of X is equal to 
\[f_2(n)=\left\{
	\begin{array}{ll}
		0  & \mbox{, if } m_j\equiv1\pmod2 \mbox{ for some } j\in\{0,1,2,\dots,n_2\} \\
		\prod_{i=1}^{n_1}(k_i+1) & \mbox{, otherwise}
	\end{array}
\right.\] 
where $n=2^{m_0} \cdot 3^{k_0} \cdot \prod_{i=1}^{n_1}p_i^{k_i} \cdot \prod_{j=1}^{n_2}q_j^{m_j}$ such that $p_i$ and $q_j$ are primes with $p_i\equiv1\pmod4$ for $i\in \{0,1,\dots,n_1\}$ and $q_j\equiv3\pmod4$ for $j\in \{0,1,\dots,n_2\}$.
\end{lemma}

\begin{proof}

Let $X$ be a map of type $[3^1,12^2;3^1,4^1,3^1,12^1]$ with associated matrix $\begin{bmatrix}a&0\\b&d\end{bmatrix}$. Now, $\psi$ sends $A$ to $B$ and $B$ to $-A$ (See Fig. \ref{Fig-1}). 
Similarly as in Lemma \ref{cm-1} we get $\psi \in \Nor(K)$ if and only if
$a\mid b$, $a\mid d$ and $ad\mid(a^2+b^2)$. Thus the polynomial is $1+x^2$.\\
Proceeding similarly as in Lemma \ref{cm-1} we get, 
$f_2(n):=\sum_{d \mid n,~n \mid d^2}\rho_{0,1}\left(\frac{d^2}{n}\right).$
and $f_2$ will be multiplicative as well. 
Now observe that,
$\rho_{0,1}(1)=1$ since $x=0$ is the solution for $1+x^2\equiv0\pmod1$.
$\rho_{0,1}(2)=1$ since $1$ is the only solution of $1+x^2=0$ in $\mathbb{Z}_2$.
$\rho_{0,1}(2^k)=0$ for all $k\in \mathbb{N}$ with $k\geq 2$. This is because when $x=2l$ for some $l\in\mathbb{N}$, $1+x^2=1+4l^2$. So, $2^k\nmid 1+x^2$ for all $k$. When $x=3l$ where $l\in\mathbb{N}\cup \{0\}$, $+x^2=1+9l^2$. When $x=3l+1$ where $l\in\mathbb{N}$, $1+x^2=9k^2+6k+2$. When $x=3l+2$ where $l\in\mathbb{N}$, $1+x^2=9k^2+12k+5$. So, $3^k\nmid(1+x^2)$ for all $x\in \mathbb{N}$. Hence, $\rho_{0,1}(3^k)=0$. Therefore,
$$f_2(2^k) =\left\{
	\begin{array}{ll}
		\sum_{l= \frac{k}{2} }^k\rho_{0,1}\left(2^{2l-k}\right) = 1  & \mbox{, if  2 $\mid$ k} \\
		\sum_{l= \frac{k+1}{2} }^k\rho_{0,1}\left(2^{2l-k}\right)=1 & \mbox{, if  2 $\nmid$ k.}
	\end{array}
          \right.
f_2(3^k)=
\left\{
	\begin{array}{ll}
		\sum_{l= \frac{k}{2} }^k\rho_{0,1}\left(2^{2l-k}\right)=1 & \mbox{, if  2 $\mid$ k} \\
		\sum_{l= \frac{k+1}{2} }^k\rho_{0,1}\left(2^{2l-k}\right)=0 & \mbox{, if  2 $\nmid$ k.}
	\end{array}
\right.$$

If $p$ is any odd prime except $3$, then, $x^2+1 \equiv 0\pmod p$ has a solution if and only if $y^2\equiv -4\pmod p$ has a solution. 
Using Legendre symbols \cite[Chap. 9]{burton2018}, $\left(\frac{-4}{p}\right)=\left(\frac{-1}{p}\right)\left(\frac{4}{p}\right)=\left(\frac{-1}{p}\right)=\left\{
	\begin{array}{ll}
		1  & \mbox{if } p\equiv1\pmod4 \\
		-1 & \mbox{if } p\equiv3\pmod4
	\end{array}
\right.$.
Hence, we have $\rho_{0,1}(p)=
\left\{
	\begin{array}{ll}
		2  & \mbox{if } p\equiv1\pmod4 \\
		0 & \mbox{if } p\equiv3\pmod4.
	\end{array}\right.$.
Putting these values in the expression of $f_2$ we get,
$ f_2(p^k)=
\left \{
	\begin{array}{lll}
	    k+1 & if p\equiv 1\pmod 4\\
		1  & \mbox{if } k\equiv 0\pmod2 ~\mbox{and}~ p\equiv 3\pmod 4 \\
		0  & \mbox{if } k\equiv 1\pmod2 ~\mbox{and}~ p\equiv 3\pmod 4.
	\end{array}
\right .$
Now, using multiplicative property of $f_2$ the assertion follows.
\end{proof}

\begin{remark}
$f_2$ is special kind of divisor function. $f_2(n)$ counts the number of divisors of $n>1$ which are congruent to 1 modulo 4.
\end{remark}

\begin{lemma}\label{cm-5}
Let $X$ be an $n$-sheeted cover such that $\mathscr{I}_X=\mathscr{I}_{E_i}$. Then number of such possible $X$ is given by $$f_6(n):=\left\{
	\begin{array}{ll}
		0  & \mbox{, if } k_i\equiv1\pmod2 \mbox{ for some } i\in\{0,1,2,\dots,n_1\}\\
		1 & \mbox{, otherwise}
	\end{array}
\right.$$ 
where 
$n=2^{m}\cdot 3^{k_0} \cdot \prod_{i=1}^{n_1}p_i^{k_i}$ such that $p_i$ is any prime other than $2$ for $i\in \{0,1,\dots,n_1\}$.
\end{lemma}
\begin{proof}

Since, $\mathscr{I}_X=\mathscr{I}_{E_2}$ so Aut($X$) will contain $\psi$ and $r_1'$. Thus, the total number of such maps is given by the number of solutions of gcd$(x^2+1,\,1-x^2)$ in $\mathbb{Z}_{\frac{d^2}{n}}$ satisfying $n\mid d^2$ for every divisor $d$ of $n$.
Let $\rho_6(n):=\#\{x\in \mathbb{Z}_n:\, x^2+1=0\,{\rm and}\,1-x^2=0\}$.
Therefore
$f_6(n):=\sum_{d \mid n,~n \mid d^2}\rho_6\left(\frac{d^2}{n}\right).$ 
By similar reason $\rho_6$ is multiplicative and as a consequence we have $f_6$ is also multiplicative. Thus it is enough to calculate for prime powers.
Now, $\rho_{6}(1)=1$ since $x=0$ is the solution for the congruences $x^2+1\equiv0\pmod1$ and $x^2-1\equiv0\pmod1$.
Since $x^2+1=x^2-1$ in $\mathbb{Z}_2[x]$, gcd$\{(x^2+1),\,(x^2-1)\}=x^2+1$.
Since $x=1$ is the only solution of $x^2+1$ in $\mathbb{Z}_2$, $\rho_{6}(2)=1$.

Since gcd$\{(x^2+1),\,(x^2-1)\}=1$ in $\mathbb{Z}_{2^k}[x]$ where $k\in \mathbb{N}\setminus\{1\}$, we have $\rho_6(2^k)=0$ for all $k$.
Since gcd$\{(x^2+1),\,(x^2-1)\}=1$ in $\mathbb{Z}_{p^k}[x]$, we have $\rho_3(p^k)=0$ for all $k\in \mathbb{N}$. Thus we get,
\\
$$f_6(2^k) =\left\{
	\begin{array}{ll}
		\sum_{l= \frac{k}{2} }^k\rho_{6}\left(2^{2l-k}\right) = 1 & \mbox{, if  2 $\mid$ k} \\
		\sum_{l= \frac{k+1}{2} }^k\rho_{6}\left(2^{2l-k}\right) = 1  & \mbox{, if  2 $\nmid$ k,}
	\end{array}
          \right. 
f_6(p^k) =\left\{
	\begin{array}{ll}
		\sum_{l= \frac{k}{2} }^k\rho_{6}\left(2^{2l-k}\right) = 1 & \mbox{, if  2 $\mid$ k} \\
		\sum_{l= \frac{k+1}{2} }^k\rho_{6}\left(2^{2l-k}\right) = 0  & \mbox{, if  2 $\nmid$ k.}
	\end{array}
          \right. $$
Now, using multiplicative property of $f_6$ the assertion follows.
\end{proof}

Hence number of $n$-sheeted $2$-uniform maps up to isomorphism there are $\frac{f_2(n)-f_6(n)}{2}+ f_6(n) = \frac{f_2(n)+ f_6(n)}{2} $.
Note that $X_0$ has 8 vertices.
Thus the number of 2-uniform maps up to isomorphism with $v$ vertices is $\Phi_{16}$  given by,
$\Phi_{16}(v) =\left\{
	\begin{array}{ll}
		0  & \mbox{, if } 8 \nmid v \\
		\frac{1}{2}[f_2 (\frac{v}{8})+f_6(\frac{v}{8})] & \mbox{, otherwise.}
	\end{array}
      \right. 
      $

\subsection{Maps of type $[4^1,8^2]$}


Let $X$ be an $n$ sheeted $2$-orbital cover of $X_0$ obtained as an orbit space of $E_{27}$. Then, $X=\frac{E_{27}}{K_{27}}$ where $K_{27}\leq H_{27}$ with $H_{27}=\langle \gamma,\delta\rangle$ and
$K_{27}=\langle w_1,\,w_2\rangle$. Since $X$ is a $2$-orbital map, $V(X)$ forms $2$-Aut$(X)$ orbits. Using Result (\ref{eqn-1}),
we can say that $V(E_{27})$ forms $2$- Nor$(K_{27})$ orbits. 
Now $V(E_i)$ forms $2$-Nor$(K_i)$  orbits when only $\tau$ and translations present in  Nor$(K_i)$. If $\mathscr{I}(X)$ contains $\tau$ and $r'_1$ then $V(E_i)$ has $2$ Nor($K_i$) orbits but instead of $r'_1$ if it contains $r_3'$ then it will become one orbital. If $\mathscr{I}(X)$ contains $\psi$ then $V(E_i)$ has $1$ Nor($K_i$) orbits. And if $\mathscr{I}(X)$ is equals to isotropy group of the tiling then also it has $1$ orbit. 
The following lemma gives the number of distinct maps having $r_3'$ in its isotropy group.
\begin{lemma}\label{lem-f4}
If $X$ be a $n$ sheeted toroidal cover of $X_0$, then the total number of maps having $r_3'$  in the automorphism group of $X$ is equal to 
\[ f_4(n)=\left\{
	\begin{array}{ll}
		\prod_{i=1}^{n_1}(k_i+1)  & \mbox{, if } k_0=0 \\
		(2k_0-1)\prod_{i=1}^{n_1}(k_i+1) & \mbox{, otherwise}
	\end{array}
\right. \]
where $n=2^{k_0}p_1^{k_1}p_2^{k_2}\dots p_{n_1}^{k_{n_1}}$ where $p_i$ is any odd prime for $i\in \{0,1,\dots,n_1\}$.
\end{lemma}

\begin{proof}
Observe that $r_3'$ sends $A$ to $B$ and $B$ to $A$. Hence, $r_1'\in \Aut(E)$. So the matrix is $\begin{bmatrix} 0&1\\1&0\end{bmatrix}$. 
Proceeding similar way as in Lemma \ref{cm-1} we get $r_3' \in Aut(X)$ if and only if $a\mid b$, $a\mid d$ and $ad\mid (b^2-a^2)$.

Therefore, the total number of required maps is given by the number of solutions of the polynomial $x^2-1$ in $\mathbb{Z}_{\frac{d^2}{n}}$ satisfying $n\mid d^2$ for every divisor $d$ of $n$. Hence, we get 
$$f_4(n):=\sum_{\substack{d \mid n,~n \mid d^2}}\rho_{0,-1}\left(\frac{d^2}{n}\right).$$
By Result \ref{result--1} we have $\rho_{0,-1}$ is multiplicative. Thus $f_4$ is so. Similarly as in Lemma \ref{cm-1} $f_4(p^k) = \sum_{l=\lceil \frac{k}{2} \rceil}^k\rho_{0,-1}\left(p^{2l-k}\right).$
Now observe that,
$\rho_{0,-1}(1)=1$, $\rho_{0,-1}(2)=1$, $\rho_3(4)=2$ and $\rho_3(2^k)=4$ for all $k\in\mathbb{N}$ with $k\geq3$.
$\rho_{0,-1}(1)=1$ since $x=0$ is the solution for $x^2-1\equiv0\pmod1$.
$\rho_{0,-1}(2)=1$ since $1$ is the only solution of $0$ in $\mathbb{Z}_2$.
$\rho_3(4)=2$ since $1$ and $3$ are the solutions of $x^2-1=0$ in $\mathbb{Z}_4$.\\
$x^2-1\equiv 0\pmod{2^k}\implies x^2\equiv1\pmod{2^k}$. We know that the congruence $x^2\equiv1\pmod{2^k}$ has exactly four incongruent solutions for $k\geq 3$. Hence, $\rho_{0,-1}(2^k)=4$. 
Proceeding in similar way as in the proof of Lemma \ref{cm-1} 
$f_4(2^k)=2k-1$ when $k\in\mathbb{N}$ and for odd prime $p$,
$f_4(p^k)=k+1$ where $k\in\mathbb{N}$. Now by multiplicative property of $f_4$ the assertion follows.
\end{proof}

By Lemma  \ref{lem-f2} we get the number of maps having $ \psi$ for  in its isotropy group is  $f_2(n)$. The number of maps having $r_1'$ in its automorphism group is given by $g(n)$ same as subsec. \ref{r_1}. since the matrix for that is $\begin{bmatrix} 1&0\\0&-1\end{bmatrix}$. The number of distinct maps having full isotropy group is $f_6(n)$. 
By Lemma \ref{lem-f4} we have the number of maps having $r_3'$ in the automorphism group of X is $f_4(n)$.

Hence number of $n$-sheeted $2$-uniform maps up to isomorphism there are $\Lambda(n) := \frac{g(n)-f_6(n)}{2}+ \frac{1}{2}[\sigma(n) - g(n) - f_4(n)-f_2(n) + 2f_6(n)]$.

Note that $X_0$ has 4 vertices.
Thus the number of 2-uniform maps up to isomorphism with $v$ vertices is $\Phi_{23}$  given by,
 $\Phi_{27}(v) = \left\{
	\begin{array}{ll}
		0  & \mbox{, if } 4 \nmid v \\
		\Lambda(n) & \mbox{, otherwise.}  
	\end{array}
      \right. $ 

\subsection{Maps of type $  [3^1, 4^1, 6^1, 4^1], [3^1,12^2], [3^1, 6^1, 3^1, 6^1]$}


Let $X$ be an $n$ sheeted $2$-orbital cover of $X_0$ obtained as an orbit space of $E_i$ where $i=21,\,22,24$ (see Fig. \ref{Fig-1}). 
Then, $X=\frac{E_i}{K_i}$ where $K_i\leq H_i$ with $H_i=\langle \gamma,\delta\rangle$ and
$K_i=\langle w_1,\,w_2\rangle$. Since $X$ is a $2$-orbital map, $V(X)$ forms $2$-Aut$(X)$ orbits. Using Result (\ref{eqn-1}),
we can say that $V(E_i)$ forms $2$- Nor$(K_i)$ orbits. 
Now $V(E_i)$ forms $3$-Nor$(K_i)$ orbits when only $\tau$ and translations present in  Nor$(K_i)$. 
If $\mathscr{I}(X)$ contains $\tau$ and $r_j$ for some $j\in \{1,2,\dots,6\}$ then $V(E_i)$ has $2$ Nor($K_i$) orbits. 
If $\mathscr{I}(X)$ contains $\sigma$ then $V(E_i)$ has $1$ Nor($K_i$) orbit. And if $\mathscr{I}(X)$ is equals to isotropy group of the tiling then also it has $1$ orbit. The following lemmas gives the number of distinct maps having $r_i$ in its isotropy group for $i=1,2,3$.

\begin{lemma}\label{lem-3}
If $X$ be a $n$ sheeted toroidal cover of $X_0$, then the total number of maps having $r_i$ in the automorphism group of $X$ is equal to $f_3(n)$ as defined in Lemma \ref{cm-1} for $i=1,2,3$.
\end{lemma}
\begin{proof}
Let the associated matrix of $X$ be $\begin{bmatrix} a & 0\\ b & d\end{bmatrix}$. 
The matrices of $r_i$ for $i=1,2,3$ are  $\begin{bmatrix} 1 & 1\\ 0 & -1\end{bmatrix}$, $\begin{bmatrix} -1 & 0\\ 1 & 1\end{bmatrix}$, $\begin{bmatrix} 0 & -1\\ -1 & 0\end{bmatrix}$ respectively.
Now proceeding in similar way as in Lemma \ref{cm-1} one can see that 
$r_1\in \Nor(K)$ if and only if {\rm (i)} $a\mid b$, {\rm (ii)} $a\mid d$ and {\rm (iii)} $ad\mid (b^2+2ab)$, 
$r_2\in$ Nor$(K)$ if and only if \rm{(i)} $a\mid b$, (ii) $a\mid d$ and (iii) $ad\mid(b^2-a^2)$, 
$r_3 \in$ Nor$(K)$ if and only if $d\mid(a+2b)$. Thus the polynomial corresponding to $r_1$ is $x^2+2x$, that for $r_2$ is $x^2-1$. For $r_3$ we need to find the number of solutions of the linear diophantine equation $dk-2b=\frac{n}{d}$ in $\mathbb{Z}_d$ for every divisor $d$ of $n$. 
From Lemma \ref{lem--3} and \ref{lem-f4} we get the number of distinct maps having $r_1$ and $r_2$ in its isotropy group is $f_3(n)$. 
Now for $r_3$ let $f_8(n):=\sum_{\substack{d \mid n}}\#\{b\in\mathbb{Z}_d: dk-2b=a,\,ad=n,\,k\in\mathbb{Z}_d\}.$
Now, observe that,
if $n=1$, then $a=1$ and $d=1$. Since $b=0$ satisfies $d\mid a+2b$, $f_8(1)=1$.
When $a=1$ and $d=2$, gcd$(d,-2)=2$ and $2\nmid a$. So, the equation $dk-2b=a$ has no solution.
When $a=2$ and $d=1$, gcd$(d,-2)=1$ and $1\mid a$. So, the equation $dk-2b=a$ has $1$ solution in $\mathbb{Z}_d$. So, $f_8(2)=1$.

If $n=2^k$, then $a=2^i$, $0\leq i\leq k$ and $d=2^{k-i}$.
When $a=1$ and $d=2^n$, gcd$(d,-2)=2$ and $2\nmid a$. So, the equation $dk-2b=a$ has no solution.
If $a=2^i$, $1\leq i\leq n-1$, gcd$(d,-2)=2$ and $2\nmid a$. So, the equation $dk-2b=a$ has $2$ solutions for each $i$.
When $a=2^n$ and $d=1$, gcd$(d,-2)=1$ and $1\mid a$. So, the equation $dk-2b=a$ has $1$ solution in $\mathbb{Z}_d$.
Thus, $f_8(2^n)=0+2+2+\dots+1=2(k-1)+1=2k-1$.

When $n=p$, we have two cases.\\
Case 1: When $a=1$ and $d=p$, gcd$(d,-2)=1$ and $1\nmid a$. So, the equation $dk-2b=a$ has $1$ solution.\\
Case 2: When $a=p$ and $d=1$, gcd$(d,-2)=1$ and $1\mid a$. So, the equation $dk-2b=a$ has $1$ solution in $\mathbb{Z}_d$.
So, $f_8(p)=2$.

If $n=p^k$, then $a=p^i$, $0\leq i\leq k$ and $d=p^{k-i}$.
For any $i$, $a=p^i$, gcd$(d,-2)=1$ and $1\nmid a$. So, the equation $dk-2b=a$ has $1$ solution each.
So, $f_8(p^k)=1+1+\dots+1=k+1$.\\
Now by multiplicative property of $f_8$ we get $f_8(n) = f_4(n)=f_3(n)$.
\end{proof}

Now, from Lemma \ref{lem-f5} we get the number of distinct maps with isotropy group equals to isotropy group of the tiling is $f_5(n)$.

Therefore, number of distinct maps having only $r_i$ in $\mathscr{I}(X)$ is $f_j(n)-f_5(n)$ for $(i,j)=(1,3),(2,4),$ $(3,8)$.
Hence, number of $n$-sheeted $2$-uniform maps up to isomorphism there are $\frac{f_3(n)+f_4(n)+f_8(n)-3f_5(n)}{3}$ $= f_3(n)-f_5(n)$.

Let there are $v_0$ vertices in $X_0$. Thus the number of 2-uniform maps up to isomorphism with $v$ vertices is $\Phi_{\ell}$  given by,
$\Phi_{\ell}(v) =\left\{
	\begin{array}{ll}
		0  & \mbox{, if } v_0 \nmid v \\
		(f_3(\frac{v}{v_0}) - f_5(\frac{v}{v_0})) & \mbox{, otherwise}  
	\end{array}
      \right. $ 
 for $(\ell,v_0) = (21,6),(22,6),(24,3)$.

\subsection{Maps of type $[3^2,4^1,3^1,4^1]$}


Let $X$ be an $n$ sheeted $2$-orbital cover of $X_0$ obtained as an orbit space of $E_{23}$. Then, $X=\frac{E_{23}}{K_{23}}$ where $K_{23}\leq H_{23}$ with $H_{23}=\langle \gamma,\delta\rangle$ and
$K_{23}=\langle w_1,\,w_2\rangle$. Since $X$ is a $2$-orbital map, $V(X)$ forms $2$-Aut$(X)$ orbits. Using Result (\ref{eqn-1}),
we can say that $V(E_{23})$ forms $2$-Nor$(K_{23})$ orbits.
Here, a map $Y$ will be $1$ orbital if and only if $\psi$ or glide reflection present in Aut$(Y)$ So whenever nither $\psi$ nor glide present Aut($Y$), $Y$ will be $2$-orbital.
Now, if both $\psi$ and glide present in Aut($X$) then we count it once in the collection of distinct $n$-sheeted maps.
If any one of $\psi$ or glide present in Aut($X$) we count it twice.
If none of them present then whenever for a map $Y$ after applying $\psi$ and glide we get equal maps we count it twice and otherwise we count it $4$ times.
Let $A(n)$ denotes number of $n$-sheeted maps up to isomorphism on which after applying $\psi$ and glide we get equal maps, and $B(n)$ denotes the number of maps up to isomorphism for which we get different maps. In both of these counting of $A(n)$ and $B(n)$ we considering those maps whose automorphism group does not contains $\psi$ and glide. From Lemma \ref{lem-f2} we get number of distinct maps having $\psi$ in its automorphism group is $f_2(n)$.  From Lemma \ref{glide} we get the the number of distinct maps having glide in its automorphism group is $g(n)$.

\begin{lemma}
Let $X$ be a $n$ sheeted toroidal map of type $[3^2,4^1,3^1,4^1]$ represented by the matrix $\begin{bmatrix}
a & 0 \\ b & d \end{bmatrix}$ having $\psi$ and glide reflection in its isotropy group. Then number of such distinct maps will be 
$$h(n):=\sum_{\substack{d \mid n,~n \mid d^2}}\rho_7\left(\frac{d^2}{n}\right) ~ where ~\rho_7(n):=\#\{x\in \mathbb{Z}_n:\, x^2+1=0\,{\rm and}\,2x=0\}.$$
\end{lemma}

\begin{proof}
$\phi$ and glide both will present in the automorphism group of a map represented by $\begin{bmatrix}
a&0\\b&d
\end{bmatrix}$ if and only if $a\mid b$, $a\mid d$, $ad\mid a^2+b^2$ and $d \mid 2b$. Hence, the assertion follows.
\end{proof}

\begin{lemma}\label{alpha}
Let $X$ be a $n$ sheeted toroidal map of type $[3^2,4^1,3^1,4^1]$ represented by the matrix $\begin{bmatrix}
a & 0 \\ b & d \end{bmatrix}$ such that after applying glide reflection and $90^{\circ}$ rotation on $X$ we get same map. Then number of distinct such maps are
$\alpha(n) = f_4(n)-h(n).$
\end{lemma}

\begin{proof}
Let $Y$ be represented by the matrix $M=\begin{bmatrix}
a&0\\b&d
\end{bmatrix}$. Let matrix of $90^{\circ}$ is denoted by $P$ and that of glide is $Q$. Then $P=\begin{bmatrix}
0&-1\\1&0
\end{bmatrix}$ and $Q=\begin{bmatrix}
-1&0\\0&1
\end{bmatrix}$. By hypothesis HNF of $PM=$ HNF of $QM$. Therefore there exists an unimodular matrix $U=\begin{bmatrix}
p&q\\r&s
\end{bmatrix}$ such that $PM=QMU$. Putting expressions of $P, Q, U$ and comparing matrix entries we get $ap=b,\, aq=d,\, bp+rd=a,\, bq+sd=0$.

If $p=0$ then $b=0,\,rd=a,\, sd=0$. Since $d\neq 0$ so $s=0$. $a$ and $d$ are both positive so $rd=a \implies r>0$. Since, $U$ is unimodular so we have $rq=\pm 1$. Therefore, $r=1$. Hence, $a=d$. Therefore, $M=\begin{bmatrix}
a&0\\0&a
\end{bmatrix}$ which contradicts the fact that $M$ does not have $90^{\circ}$ rotation in its automorphism group. So, $p\neq 0$. Since $a\neq 0$ so $b\neq 0$. $bp+rd=a \implies 1-p^2=rq$ and $bq+sd=0 \implies p(p+s)=0 \implies s=-p$. For a fixed $d$  if we know $p$ then we can derive values of other unknowns. So number of such possible maps is given by $\sum_{d\mid n,~n \mid d^2}[{\rm Number~ of ~solutions~ of} x^2-1\equiv 0 \pmod{\frac{d^2}{n}}]$(putting $q=d^2/n$ as $a=n/d)$. We will see in Lemma \ref{lem-f4}, this quantity is equals to $f_4(n)$. Since we are considering those maps which does not have $90^{\circ}$ and glide in its automorphism group so to get required maps we have to subtract $h(n)$ from above. 
Hence, $\alpha(n) = f_4(n)-h(n).$ 
\end{proof}

Now we have,
$$2A(n) + 4B(n) = \sigma(n) - f_2(n) - g(n) + h(n)~~~~~~~~(*)$$
By Lemma \ref{alpha} we have $2A(n) = \alpha(n) \implies A(n) = \frac{\alpha(n)}{2}$.
Putting the value of $A(n)$ in $(*)$ we get,
$B(n) = \frac{1}{4} [\sigma(n) - f_2(n) - g(n) + h(n) - \alpha(n) ] $.
Thus number of $2$-uniform maps up to isomorphism is, 
$$A(n) + B(n) = \frac{1}{4}[\sigma(n) - g(n) -f_2(n) + h(n) + \alpha(n)] =\frac{1}{4}[\sigma(n) - g(n) -f_2(n) + f_3(n)] .$$

Note that $X_0$ has 4 vertices.
Thus the number of 2-uniform maps up to isomorphism with $v$ vertices is $\Phi_{23}$  given by,
 $\Phi_{23}(v) = \left\{
	\begin{array}{ll}
		0  & \mbox{, if } 4 \nmid v \\
		\frac{1}{4}[\sigma (\frac{v}{4})-g(\frac{v}{4}) - f_2(\frac{v}{4})+f_3(\frac{v}{4})] & \mbox{, otherwise.}  
	\end{array}
      \right.$ 

\subsection{Maps of type $[3^4,6^1]$}


Let $X$ be an $n$ sheeted $2$-orbital cover of $X_0$ obtained as an orbit space of $E_{25}$. Then, $X=\frac{E_{25}}{K_{25}}$ where $K_{25}\leq H_{25}$ with $H_{25}=\langle \gamma,\delta\rangle$ and
$K_{25}=\langle w_1,\,w_2\rangle$. Since $X$ is a $2$-orbital map, $V(X)$ forms $2$-Aut$(X)$ orbits. Using Result (\ref{eqn-1}),
we can say that $V(E_{25})$ forms $2$- Nor$(K_{25})$ orbits. It can be observed that from the figure that $\gamma,\,\delta,\,\tau,\sigma \in \Aut(E_{25})$. Then, following the calculations done for the tiling $E_1$, 
we can say that $G_1=\langle \gamma,\,\delta,\,\tau\rangle\leq \Nor(K_{25})$. It is easy to see that $V(E_{25})$ forms $3$-$G_1$ orbits.
Let $G_2=\langle \gamma,\,\delta,\,\tau, \,\sigma\rangle$. If $G_{25}\leq \Nor(K_{25})$ and we observe that $V(E_{25})$ forms $1$-$G_2$ orbit. Since there are no other symmetries in Aut$(E_{25})$, the total number of $2$-orbital $n$-sheeted maps is given by $0$.

Thus,
 $\Phi_{25}(v) = 0$ for all $v\in \mathbb{N}$. 

\subsection{Maps of type $[4^1,6^1,12^1]$}\label{last}


Let $X$ be an $n$ sheeted $2$-orbital cover of $X_0$ obtained as an orbit space of $E_{26}$. Then, $X=\frac{E_{26}}{K_{26}}$ where $K_{26}\leq H_{26}$ with $H_{26}=\langle \gamma,\delta\rangle$ and
$K_{26}=\langle w_1,\,w_2\rangle$. Since $X$ is a $2$-orbital map, $V(X)$ forms $2$-Aut$(X)$ orbits. Using Result (\ref{eqn-1}),
we can say that $V(E_{26})$ forms $2$- Nor$(K_{26})$ orbits.
Now $V(E_i)$ forms $6$-Nor$(K_i)$  orbits when only $\tau$ and translations present in  Nor$(K_i)$. If $\mathscr{I}(X)$ contains $\tau$ and $r_j$ for some $j\in \{1,2,\dots,6\}$ then $V(E_i)$ has $3$ Nor($K_i$) orbits. If $\mathscr{I}(X)$ contains $\rho$ then $V(E_i)$ has $2$ Nor($K_i$) orbits. And if $\mathscr{I}(X)$ is equals to isotropy group of the tiling then also it has $1$ orbits.
By Lemma \ref{cm-1} number of distinct maps having $\rho$ in its automorphism group is $f_1(n)$ and by Lemma \ref{lem-f5} number of distinct maps having isotropy group equals to isotropy group of tiling is $f_5(n)$.
Hence number of $n$-sheeted $2$-uniform maps up to isomorphism there are $\frac{f_1(n)-f_5(n)}{2} $.

Note that $X_0$ has 12 vertices.
Thus the number of 2-uniform maps up to isomorphism with $v$ vertices is $\Phi_{23}$  given by,
 $\Phi_{26}(v) =\left\{
	\begin{array}{ll}
		0  & \mbox{, if } 12 \nmid v \\
		\frac{1}{2}[f_1(\frac{v}{12}) - f_5(\frac{v}{12})] & \mbox{, otherwise.}
	\end{array}
      \right. $ 

\begin{proof}[Proof of theorem \ref{thm-1} {\rm (a,b)}]
The proof follows from subsec. \ref{first} to \ref{last}. The explicit functions $\Phi_{ell}$ are given in respective subsections. This completes the proof of part (a) of Theorem \ref{thm-1}.

Now, to get explicit maps we do the following. Suppose number of vertices be $v$ and we want to find 2-uniform maps of type $\ell$. Let $v_0$ be the number of vertices of minimal map of that type. Then if $v_0 \nmid v$, $\Phi_{ell}(v) = 0$. So, there is no such map. Otherwise, let $v=nv_0$. Then list out all HNF of $2 \times 2$ integer matrices with determinant $n$. Note also that which symmetries of the tiling gives 2 vertex orbits.
Then according to the type use corresponding Lemmas in Sec. \ref{proof} to find the maps upto isomorphism. This completes the proof of part (b) of Theorem \ref{thm-1}.
\end{proof}

We end this section with some examples.

\begin{example}
For maps of type $[3^2,4^1,3^1,4^1]$ number of 2-uniform maps with 32 vertices is 3. The HNF of those maps will be $\begin{bmatrix}2&0\\1&4\end{bmatrix}, \begin{bmatrix}1&0\\3&8\end{bmatrix}$ and $\begin{bmatrix}1&0\\1&8\end{bmatrix}$. 
\end{example}

\begin{example}
For the same type number of 2-uniform maps with 36 vertices is 3. The HNF of those maps will be $\begin{bmatrix}1&0\\1&9\end{bmatrix}, \begin{bmatrix}1&0\\2&9\end{bmatrix}$ and $\begin{bmatrix}1&0\\3&9\end{bmatrix}$. 
\end{example}

\section{Asymptotic behaviour}\label{ab}

In this section we discuss about the asymptotic behaviour of $\Phi_{\ell}(v)$ as $v\to \infty$. We will give an asymptotic upper bound of $\Phi_{\ell}(v)$. 
\begin{proof}[Proof of theorem \ref{thm-1} {\rm (c)}]
First recall the Gronwall's theorem \cite{gronwall} stated as follows:
\begin{equation}
    \Lim \frac{\sigma(v)}{v \ln \ln v} = e^{\gamma}
\end{equation}
where $\gamma$ is Euler-Mascheroni constant. It can be easily seen that $\Lim \frac{f_k(v)}{v \ln \ln v} = 0$ for all $k$.

Therefore we get, 
\begin{equation}\label{lim}
\Lim \frac{\Phi_{27}(v)}{v \ln\ln \frac{v}{4}} = \frac{1}{8} e^{\gamma}.
\end{equation}
Thus $\Phi_{27}(v)$ is bounded above by $\frac{1}{8} e^{\gamma} v \ln\ln \frac{v}{4}$ for large enough $v$.
Following remark follows from \ref{lim}.
\begin{remark}
There exists $v\in \mathbb{N}$ such that $\Phi_{27}(v) > v$. More precisely: The number of 2-uniform toroidal maps up to isomorphism of type $[4^1,8^2]$ with $v$ vertices can be grater than $v$. 
\end{remark}
The exactly similar remark can be drawn for $\Phi_{\ell}$ when $\ell = 3,4,8,12,13,15,23$ also.

From \cite{gronwall} we get, 
\begin{equation}
\Lim \frac{\tau(v)}{\frac{\ln 2 \ln v}{\ln \ln v}}=1.
\end{equation}
Now note that $f_i(v) \le \tau(v)$ for all $i=1,2,\dots,6$ and for all $v\in \mathbb{N}$. Therefore we have,
$$\Phi_{\ell}(v) \le \tau(\frac{v}{v_0}) \le \tau(v) \implies \ln (\Phi_{\ell}(v)) \le \ln (\tau(v)) \implies \Lim \frac{\ln (\Phi_{\ell}(v))}{\frac{\ln 2 \ln v}{\ln \ln v}} \le 1$$
for $\ell = 1, 2, 5, 6, 7, 11, 14, 16, 17, 19, 20,21, 22, 24, 26$.
Thus $\ln (\Phi_{\ell}(v))$ is bounded above by $\frac{\ln 2 \ln v}{\ln \ln v}$ for large enough $v$. Hence $\Phi_{\ell}(v)$ is bounded above by $e^{\frac{\ln 2 \ln v}{\ln \ln v}}$ for large enough $v$.

Now, from Remark \ref{remark1} $g(n) \le 2\tau(n)$. Thus $\Phi_{\ell}(v)$ is bounded above by $2\frac{\ln 2 \ln v}{\ln \ln v}$ for large enough $v$ when $\ell = 9,10,18$.

For all $\ell$ the lower bound of $\Phi_{\ell}$ is the constant zero function.
\end{proof}

\section{Acknowledgements}
The second author is partially supported by Science and Engineering Research Board (SERB), DST, India (SRG/2021/000055) and both authors are supported by NBHM, DAE (No. 02011/9/2021-NBHM(R.P.)/R$\&$D-II/9101).

{\small
\bibliographystyle{plain}
\bibliography{article_V5}

}

\end{document}